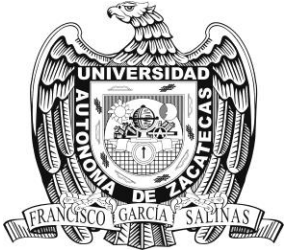 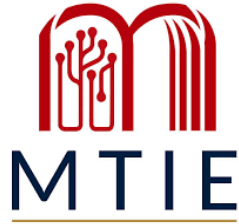 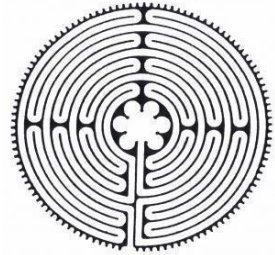

# UNIVERSIDAD AUTÓNOMA DE ZACATECAS "FRANCISCO GARCÍA SALINAS"

# UNIDAD ACADÉMICA DE DOCENCIA SUPERIOR MAESTRÍA EN TECNOLOGÍA INFORMÁTICA EDUCATIVA

"Material didáctico de Lógica Proposicional para Estructuras Discretas"

Presenta:
Margarita Carrera Fournier

Directora:
Dra. Verónica Torres Cosío

Zacatecas, Zac., julio 2022

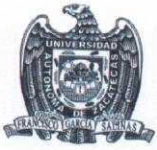 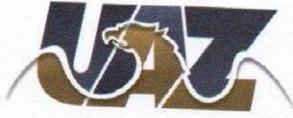 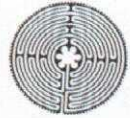 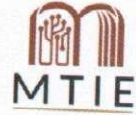

**Asunto:** Autorización de Impresión de Trabajo
No. Oficio MTIE 031/2022

**C. CARRERA FOURNIER MARGARITA**
Candidato (a) a Grado de Maestría en
Humanidades y Procesos Educativos
P R E S E N T E

Por este conducto, me permito comunicar a usted, que se le autoriza para llevar a cabo la impresión de su trabajo de tesis:

*"Material Didáctico de Lógica Proposicional para Estructuras Discretas".*

Que presenta para obtener el Grado de Maestría.

También se le comunica que deberá entregar a este Programa Académico (2) dos copias de su tesis a la brevedad posible.

Sin otro particular de momento, me es grato enviarle un cordial saludo.

**A T E N T A M E N T E**
Zacatecas, Zac., a 18 de agosto del 2022

**Dra. Verónica Torres Cosío**
**Responsable del Programa de la MTIE**

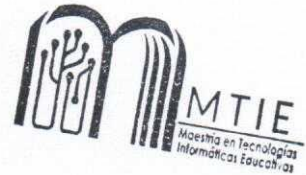

c.c.p.- Alumno
c.c.p.- Archivo

Dra. Verónica Torres Cosío
Responsable de la MTIE
P R E S E N T E

En respuesta al nombramiento que me fue suscrito como directora de tesis de la alumna: **Margarita Carrera Fournier** cuyo título de su trabajo se enuncia: **"Material Didáctico de Lógica Proposicional para Estructuras Discretas".**

**Hago constar que ha cubierto los requisitos de dirección y corrección satisfactoriamente**, por lo que está en posibilidades de pasar a la disertación de su trabajo de investigación para certificar su grado de Maestro (a) en Tecnología Informática Educativa. De la misma manera no existe inconveniente alguno para que el trabajo sea autorizado para su impresión y continué con los trámites que rigen en nuestra institución.

Se extiende la presente para los usos legales inherentes al proceso de obtención del grado del interesado.

A T E N T A M E N T E
Zacatecas, Zac., a 16 de agosto del 2022

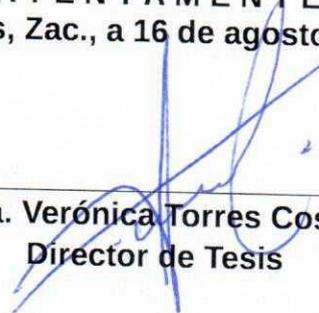

Dra. Verónica Torres Cosío
Director de Tesis

c.c.p.- Interesado
c.c.p.- Archivo

**Dedicatoria**

Dedico este trabajo al motor de mi vida, mis hijos André y Rodrigo, siempre serán mi inspiración y motivación para crecer profesionalmente y como persona en todos los sentidos.

A mi esposo por ser quién me motivo a emprender este nuevo reto, pero sobre todo por su apoyo, motivación, amor y no dejarme rendir. Definitivamente sin él no lo habría logrado.

A mis padres porque desde pequeña han creído en mí y me motivan a seguir mis sueños.

A mis hermanas, cuñados y sobrinos por la confianza y palabras de aliento.

Gracias familia, los amo.

# Resumen

Uno de los temas que presenta mayor dificultad en la asignatura de Matemáticas Discretas es el tema de Lógica proposicional, por ello el presente trabajo tuvo como objetivo facilitar el aprendizaje de la Lógica Proposicional mediante la implementación de material didáctico con el uso de tecnología educativa que contribuya en el logro de los objetivos establecidos en el temario de la asignatura de Matemáticas Discretas, debido a que los conocimientos que aporta en los Ingenieros en Computación son fundamentales para su desarrollo laboral. El estudio se consideró del tipo cuasiexperimental de corte descriptivo, en donde se utilizó un instrumento de diagnóstico y otro de satisfacción. A través de la descripción se presenta tanto el rediseño de material didáctico con base al modelo de diseño instruccional ASSURE como su implementación mediante el proceso de intervención educativa. Se contó con la participación de dos grupos de estudiantes no equivalentes de la Facultad de Ingeniería de la Universidad Nacional Autónoma de México, quienes fueron elegidos por conveniencia. Los resultados revelaron que se logró avanzar y cumplir en un alto porcentaje del tema de lógica proposicional, esto a pesar de haberse presentado otros factores no contemplados, como el cambio de modalidad de presencial a en línea, derivado de la pandemia por COVID-19.

***Palabras clave:*** lógica proposicional, material didáctico, intervención educativa, TIC

# Índice

## Índice de figuras

**Índice de tablas**

# Capítulo I. Introducción

La lógica proposicional es un instrumento que puede ser empleado en diferentes campos de interés, en particular para las carreras con un énfasis informático aporta la habilidad de desarrollar el razonamiento matemático así mismo ayuda a obtener una mejor comprensión de la realidad representándolo mediante proposiciones lógicas. Por lo cual el propósito de este trabajo de investigación es facilitar su aprendizaje en particular de los estudiantes de la carrera de Ingeniería en Computación de la Universidad Nacional Autónoma de México.

El presente trabajo da inicio con los antecedentes o estudios realizados previamente similares al que se desarrolló, posteriormente se describe el marco contextual donde se desarrolló, así como el planteamiento de la problemática detectada, los objetivos, la justificación, los alcances y limitaciones.

El segundo capítulo es el marco teórico que está constituido por tres capítulos principales, Lógica proposicional, Métodos lógicos para la demostración de argumentos y Material didáctico mediado por las TIC.

En el siguiente capítulo se plantea la metodología donde se detallan los pasos que se siguieron para desarrollar el presente proyecto, desde la definición de la metodología, la población, las técnicas e instrumentos utilizados hasta el diseño que se implementó.

En el penúltimo capítulo se detallan los resultados obtenidos que van en congruencia con el objetivo general concretamente es facilitar el aprendizaje de la lógica proposicional mediante la implementación de material didáctico con el uso de tecnología educativa asimismo los objetivos específicos que sustentan al objetivo principal.

Por último, se describen las conclusiones logradas del proyecto donde a partir de los resultados se puede decir que dicha una estrategia empleada es una buena herramienta complementaria para beneficiar el aprendizaje de los estudiantes.

## 1.1 Antecedentes

Realizando una investigación sobre proyectos similares en el plano internacional existen trabajos realizados en Cuba y Colombia, pero enfocados al uso de redes sociales como Facebook, Twitter o investigaciones sobre el impacto de la manera en que los docentes enseñan las Matemáticas. Lo más parecido a lo que se desarrolló en este trabajo es el de la Universidad Nacional de Costa Rica, donde Vilchez (2013) realizó un trabajo llamado Cuadernos interactivos para un curso de Estructuras Discretas.

Su desarrollo consistió en elaborar todos los cuadernos interactivos para el abordar los contenidos vinculados a un curso de Estructuras Discretas para Informática, que se imparte a los estudiantes de la carrera Ingeniería en Sistemas de Información de la Universidad Nacional de Costa Rica y Costa Rica. Entendiendo como cuaderno interactivo a una aplicación informática que le permite al estudiante profundizar cada uno de los temas del curso, utilizando como principales recursos animaciones y videos educativos

Su objetivo fue analizar el impacto de la red social Facebook como un entorno educativo, se planteó la necesidad de diseñar una serie de materiales interactivos para acompañar el aprendizaje de los estudiantes durante una experiencia de implementación.

Desarrollaron 8 cuadernos interactivos destinados a la explicación en video del curso de Estructuras Discretas apoyándose de la compañía Livescribe, quien les proporciono la tecnología y soporte on-line a sus usuarios con la finalidad de publicar en la web un formato de archivo denominado pencast con una conectividad directa a perfiles creados en Facebook. Además, emplearon animaciones en Flash llamadas "quices".

Gracias a que realizaron una encuesta de percepción por parte de los estudiantes y de un taller de valoración de estos materiales con profesores de la cátedra con la intención de determinar su calidad académica y pedagógica, el resultado fue exitoso, aunque demandante.

Siguiendo en el plano nacional en la Universidad de Guadalajara realizaron un trabajo llamado Objetos de aprendizaje como recursos didácticos para la enseñanza de matemáticas, donde Aragón et al. (2009) buscaban facilitar la enseñanza de conceptos matemáticos mediante el uso de un objeto de aprendizaje apoyado en recursos tecnológicos. Su muestra fueron 170 estudiantes de seis grupos de nivel licenciatura de diferentes cursos y áreas disciplinares. Las instituciones participantes fueron la Escuela Normal Básica Miguel Hidalgo; Benemérita y Centenaria Escuela Normal del Estado de Chihuahua; ITESM, campus Monterrey; Universidad de Quintana Roo, campus Cozumel; y Escuela de Administración San Pedro de la Universidad Autónoma de Coahuila. Una característica destacable fue que todos cursaron al mismo tiempo la asignatura y aunque el resultado de esta estrategia fue innovadora y positiva fue aplicada para las matemáticas en general no la rama de las Matemáticas Discretas.

En cuanto a la Universidad Nacional Autónoma de México, no hay trabajos de investigación similares registrados. Generalmente están enfocados a las matemáticas en general, pero no en específico para Matemáticas Discretas.

## 1.2 Marco contextual

La Facultad de Ingeniería de la Universidad Autónoma de México es la sede donde se llevará acabo la intervención y se encuentra dividida en la sede central ubicada en el Circuito Escolar en Ciudad Universitaria, y la División de Ciencias Básicas, ubicada al sur de las instalaciones deportivas de la UNAM.

En la Facultad de Ingeniería de la Universidad Nacional Autónoma de México se imparten 15 carreras a nivel licenciatura dentro de las cuales está la de Ingeniería en Computación. Cuenta con 2,171 académicos (el 72% de los académicos son hombres y el 28% son mujeres) para impartir las asignaturas de las diferentes carreras y específicamente 547 para la carrera de Ingeniería en Computación.

A la Facultad de Ingeniería de la Universidad Nacional Autónoma de México aproximadamente ingresan 450 estudiantes en cada generación de Ingenieros en Computación, donde el 80% son hombres y el 20% mujeres.

Al ingresar los estudiantes necesitan tener conocimientos de matemáticas en álgebra, geometría analítica y cálculo diferencial e integral de funciones de una variable, física, particularmente en lo que respecta a temas relacionados con mecánica clásica, así como conocimientos generales de química y de computación. Es también conveniente que posea conocimientos de inglés, por lo menos a nivel de comprensión de textos. Por lo que respecta a las habilidades, es importante que tenga disposición para el trabajo en equipo, capacidad de análisis y síntesis, y de adaptación a situaciones nuevas, así como espíritu creativo.

En general se busca que los egresados de la Facultad de Ingeniería posean: capacidades para la innovación, potencial para aportar a la creación de tecnologías y actitud emprendedora, con sensibilidad social y ética profesional; y con potencialidad y vocación para constituirse en factor de cambio.

*Misión*

Generar recursos humanos en ingeniería con una formación integral de excelencia académica, con un sentido ecológico, ético y humanista que los compromete a mantenerse actualizados permanentemente, capaces de resolver problemas de forma creativa e innovadora en el ámbito de su competencia, así como de realizar investigación científica y aplicada acorde a las necesidades de la sociedad y de impacto en el desarrollo nacional.

*Visión*

La Facultad de Ingeniería es una institución educativa de excelencia, referente nacional y de prestigio internacional. Formadora de profesionales, en los niveles de licenciatura y posgrado, altamente competitivos y demandados por los sectores productivos debido al dominio de sus conocimientos en ingeniería. Su personal académico es líder en su campo, con una alta productividad científica y tecnológica, tal que le permite realizar investigación de punta para resolver los problemas nacionales.

Para lograr esta visión, el personal académico debe fomentar la participación de los estudiantes en proyectos de investigación y publicar sus avances en materia de generación de nuevo conocimiento en revistas arbitradas nacionales e internacionales. Además, la Facultad de Ingeniería tiene que establecer estrategias de largo plazo, basadas en la ética, el trabajo colaborativo, la honestidad, la perseverancia, la equidad, la responsabilidad y la racionalidad en el uso de los recursos que le permitan alcanzar su visión y consolidarse en ella en un ambiente académico-administrativo de primer mundo.

## 1.3. Planteamiento del problema

**Análisis del problema**

Durante 14 años a través de la observación a estudiantes que han cursado la asignatura de Matemáticas Discretas de manera presencial, se identificó que presentaban mayor problema en el tema de lógica proposicional sin importar los antecedentes teóricos iniciales de los estudiantes o si se incrementa el número de ejercicios realizados durante la clase, así como los elaborados como parte de sus tareas, el tema más complicado para los estudiantes es el de lógica proposicional.

En cuanto a las calificaciones de la evaluación del tema de lógica proposicional se identificó mediante un examen parcial la dificultad en la comprensión y resolución de ejercicios a comparación de los otros temas que contempla el temario de la asignatura de Matemáticas Discretas, no obstante que en los últimos cuatro semestres se dividió el tema y en lugar de una sola evaluación, se realizó en dos.

En este sentido también se les ha invitado a que de manera voluntaria realicen ejercicios fuera de clase y su desarrollo sea validado y regresado con las observaciones pertinentes si fuera el caso, sin embargo, los estudiantes refieren que necesitan más ejercicios con solución para que de manera independiente puedan resolverlos por ellos mismos.

Como consecuencia de lo expuesto anteriormente es que se busca la implementación de material didáctico con el uso de tecnología educativa para lograr los objetivos establecidos en el temario de la asignatura.

**Planteamiento del problema**

El egresado de Ingeniería en Computación tiene que contar con sólidas bases científicas y fundamentos tecnológicos, que le permitan comprender, analizar, diseñar, organizar, producir, operar y dar soluciones prácticas a problemas relacionados con las áreas de Organización de Sistemas Computacionales, Ingeniería en Software y Tecnologías de Información. Adicionalmente con base en el campo de profundización seleccionado, tendrá conocimientos en algunas áreas tecnológicas tales como: procesamiento digital de datos y control de procesos, sistemas de programación tanto de base como de aplicación, desarrollo e investigación en las ciencias de la computación, sistemas de comunicación y seguridad tanto informática como de redes de datos, sistemas de bases de datos, sistemas inteligentes, y sistemas de cómputo gráfico, entre otras.

A pesar de que la lógica proposicional es de los temas con mayor grado de dificultad de aprendizaje de la materia de Estructuras Discretas, los estudiantes de la carrera de Ingeniería en Computación de la Facultad de Ingeniería de la Universidad Nacional Autónoma de México no disponen de material, ni literatura suficiente con ejercicios resueltos con un mayor grado de dificultad.

Cada semestre se lleva a cabo juntas de academia (de acuerdo con la asignatura), aproximadamente tres durante el semestre, en las cuales se exponen el cumplimiento del temario de la asignatura, el estado que guardan los grupos, presentación de sugerencias didácticas a desarrollar durante el semestre, entre otros temas.

De los puntos mencionados en estas reuniones se tiene:

- Ninguno de los integrantes de la academia terminaba el temario.

- Deficiencias de los antecedentes académicos de los estudiantes.
- Apatía de los estudiantes.

Derivado de lo anterior, se realizó un cambio de plan de estudios, el cual no sólo se modificó en cuanto a contenido, si no también se redujeron las horas. Este plan de estudios se puso en marcha desde el 2017. Esta solución no resolvió el que se concluyera con el temario, además de que al reducir las horas destinadas para impartir la materia provocó que no se tenga el tiempo suficiente para dar los ejercicios necesarios que refuercen el conocimiento de los estudiantes.

## 1.4 Objetivos

**Objetivo general**

Facilitar el aprendizaje de la lógica proposicional mediante la implementación de material didáctico con el uso de tecnología educativa que contribuya en el logro de los objetivos establecidos en el temario de la asignatura.

**Objetivos específicos**

1. Rediseñar ejercicios de lógica proposicional por diferentes métodos de solución.
2. Producir materiales didácticos por medio de herramientas de software educativas, en apego a los objetivos de la asignatura.
3. Implementar los materiales didácticos en apego a los objetivos de la asignatura.
4. Analizar el rendimiento académico después de haber implementado el material didáctico.
5. Medir el nivel de satisfacción del curso mediado por tecnología.

## 1.5 Preguntas de investigación

**Pregunta general**

¿La implementación del material didáctico sobre el tema de lógica proposicional de la materia de Estructuras Discretas facilita su aprendizaje de acuerdo con los objetivos establecidos en el temario de la asignatura?

**Preguntas específicas**

1. ¿El rediseño de ejercicios de lógica proposicional por diferentes métodos de solución de forma detallada y desglosando paso a paso su solución se apegó a los objetivos de la asignatura?
2. ¿La producción de materiales didácticos por medio de herramientas de software educativo, se apegó a los objetivos de la asignatura?
3. ¿La implementación de los materiales didácticos en apegó a los objetivos de la asignatura?
4. ¿El análisis del rendimiento académico de acuerdo con los objetivos de la asignatura reflejó la eficacia de la elaboración del material didáctico?
5. ¿La medición del curso mediado por tecnología fue satisfactoria?

## 1.6 Justificación

Es importante el impacto que tiene la asignatura de Estructuras Discretas porque da soporte a las habilidades y conocimientos necesarios a los futuros ingenieros en Computación, considerando que puede asumir varios roles, desempeñarse con grupos multidisciplinarios así como con clientes en diferentes áreas, debe tener la capacidad de razonar, resolver situaciones complejas en el área de trabajo en donde se desenvuelva, además las soluciones planteadas deben ser estratégicas conforme a las necesidades del mercado actual y ser adaptables a los avances tecnológicos.

El aprendizaje de las matemáticas en general es considerado difícil, esto se atribuye a factores como la estrategia de enseñanza, el requerimiento del razonamiento deductivo por parte del estudiante, los diferentes estilos de aprendizaje, entre otros. La principal finalidad de la elaboración del material didáctico guiado con el uso de la tecnología educativa es facilitar la construcción del aprendizaje significativo del estudiante de la lógica proposicional, de acuerdo con los objetivos establecidos en el temario de la materia. Una de las ventajas de esta nueva estrategia es la disponibilidad del material porque van a poder consultarlo y revisarlo las veces que consideren necesario para comprender el tema a su ritmo y de forma independiente, además de no interferir en el tiempo establecido de la clase. También se busca favorecer que se abarquen todos los temas del temario, ya que se verían menos ejercicios en clase.

Metodológicamente hablando se busca beneficiar inicialmente a los estudiantes del grupo 01 de Estructuras Discretas, pero puede ser implementada por todos los profesores de la Facultad de Ingeniería de la UNAM u otras universidades donde se imparta la materia buscando disminuir los índices de reprobación. Inclusive esta nueva estrategia podría ayudar a ampliar la cobertura del modelo tradicional donde generalmente se atienden a 30 estudiantes por grupo consiguiendo hasta 50 estudiantes por grupo.

## 1.7 Alcances y limitaciones

**Alcances**

El alcance del proyecto fue en un principio para los estudiantes de la materia de Estructuras Discretas del grupo 01, durante el periodo comprendido del 21 de septiembre del 2020 al 24 de julio del 2021; sin embargo, la finalidad es que lo pueda emplear cualquier estudiante de carreras afines a sistemas computacionales que incluya en su plan de estudios la asignatura de Estructuras Discretas (también conocida como Matemáticas Discretas).

**Limitaciones**

*No existe colaboración entre profesores.*

Cada semestre se realizan juntas académicas en las que se plantean varios puntos relacionados con la materia de Estructuras Discretas, dentro de esos puntos se invita a colaborar con material para ponerlo a disposición de los estudiantes, sin embargo, ninguno de los profesores tiene dicha disposición.

Los conocimientos de los estudiantes no son uniformes (a pesar de seguir un temario), varía entre profesores por lo que este depende con quien tomen la materia.

Según Becerra, "Comprender el valor de la colaboración como competencia esencial dentro de la sociedad del conocimiento ayuda más no garantiza el éxito, es conveniente trazar líneas de acción en cada institución educativa donde se plasmen acciones, estrategias e indicadores para vivir la colaboración entre docentes". (s.f.)

*Rezago tecnológico de los docentes.*

Aunque la Facultad cuenta con la infraestructura tecnológica como plataforma LMS, no todos hacen uso de ella debido a que no saben cómo usarla en la impartición de su asignatura. Por lo que la tecnología es subutilizada y se pierden de la posibilidad de obtener y/o mejorar sus resultados en cuanto al aprendizaje de los estudiantes y sus prácticas docentes.

El docente debe contar con una capacitación y actualización continua. Según Echeverría (2020), "Hace falta de una actualización informativa para estar a la par de lo que se vive día a día, es una desventaja académica que el maestro no domine medios informáticos".

*Escaso material didáctico y el que existe tiene deficiencias.*

Existe un libro elaborado por el coordinador de la asignatura que abarca solo un capítulo de la asignatura con problemas básicos.

Los estudiantes buscan material en otra literatura o en internet en sitios donde la información no siempre es la adecuada.

Contar con conocimientos para elaborar material didáctico no sólo en cuanto a contenido sino a la estructura formal del mismo.

*Se usa las TIC, pero el material es excesivamente básico.*

Existe una plataforma donde se están realizando videos explicativos de algunos temas, pero los ejemplos son demasiado básicos y solo lo saben los estudiantes de la profesora que lleva el proyecto.

Al ser problemas muy básicos los pocos estudiantes que saben de dicho sitio pierden el interés y dejan de usarla.

El uso adecuado de las TIC dará mejor resultado como bien menciona Cabero (2007), "son solamente medios y recursos didácticos, que deben ser movilizados por el profesor cuando les puedan resolver un problema comunicativo o le ayuden a crear un entorno diferente y propicio para el aprendizaje".

# Capítulo 2. Fundamento teórico

El presente capítulo contiene el fundamento teórico que como bien considera Hernández (2008) es una de las fases más importantes de un trabajo de investigación, la cual consiste en desarrollar la teoría que va a fundamentar el proyecto con base al planteamiento del problema que se ha realizado.

Así pues, esto le dará elementos al lector para que comprenda la temática que se desarrolla en el presente trabajo abordando desde la definición, aplicación, la relación con otros ejes del conocimiento hasta algunos métodos de demostración que existen para demostrar su validez.

Por otro lado, se explicará brevemente los beneficios del material didáctico mediado por las TIC, el uso del con el modelo instruccional ASSURE aplicado a la educación a distancia, así como la importancia del aprendizaje significativo.

## 2.1 Lógica proposicional

En principio cabe retomar la definición de Matemáticas Discretas, empleando las palabras de Rosen (2004) es la disciplina que se encarga, del estudio de los objetos discretos, si se entiende por discretos a los elementos distintos o inconexos. Son la base de todo lo relacionado con los números naturales o conjuntos numerables, por tal motivo resulta ser fundamental para las ciencias de la computación, porque los datos que se almacenan en la computadora son unos y ceros, es decir, de manera discreta. Es importante esclarecer en qué consiste, debido a que el tema de la lógica proposicional forma parte de los contenidos que se abordan en la materia de Matemáticas Discretas.

*2.1.1 ¿Qué es la lógica proposicional?*

Como plantea Arnaz (2019), la lógica proposicional es la parte de la lógica que estudia las formas en que se relacionan unas proposiciones con otras y, sobre todo, la relación que se da entre las proposiciones que componen un razonamiento.

Es importante destacar que la lógica proposicional es uno de los temas de las Matemáticas Discretas por ser la base del razonamiento matemático; ya que posee reglas que especifican el significado de los enunciados matemáticos y que se utilizan para distinguir entre argumentos válidos y no válidos.

*2.1.2 Aplicaciones y su relación con otros ejes del conocimiento.*

La lógica proposicional tiene aplicaciones en áreas de ciencias de la computación en prácticas de diseño de equipos informáticos, la especificación de sistemas, la inteligencia artificial y en los lenguajes de programación, ya que, si hablamos de la solución o soluciones de un problema específicamente computacional, debe ser resuelto mediante algoritmos y estos deben ser eficientes para poder ser ejecutados con la capacidad de cálculo necesaria. Desde el punto de vista de Farré et al. (2012) afirman que:

> Los informáticos necesitan analizar las propiedades lógicas de sus sistemas mientras los diseñan, desarrollan, verifican y mantienen, especialmente cuando se trata de sistemas críticos económicamente (un fallo tendría un alto coste económico), o críticos en seguridad o privacidad. Pero también en sistemas cuya eficiencia resulta crítica, los análisis lógicos pueden ser reveladores y, en todos los tipos de sistemas, los métodos y herramientas basados en la lógica pueden mejorar la calidad y reducir costes. (p. IX)

En relación con lo anterior se puede agregar que la lógica proposicional nos ayudar a disminuir la cantidad de errores porque nos enseña a construir e interpretar sentencias sobre el mundo que nos rodea en un sentido lógico en base a nuestro raciocinio, en este sentido es por ello su relación con otras ciencias como lo muestra la Figura 1.

**Figura 1**

*Relación de la lógica con otras ciencias.*

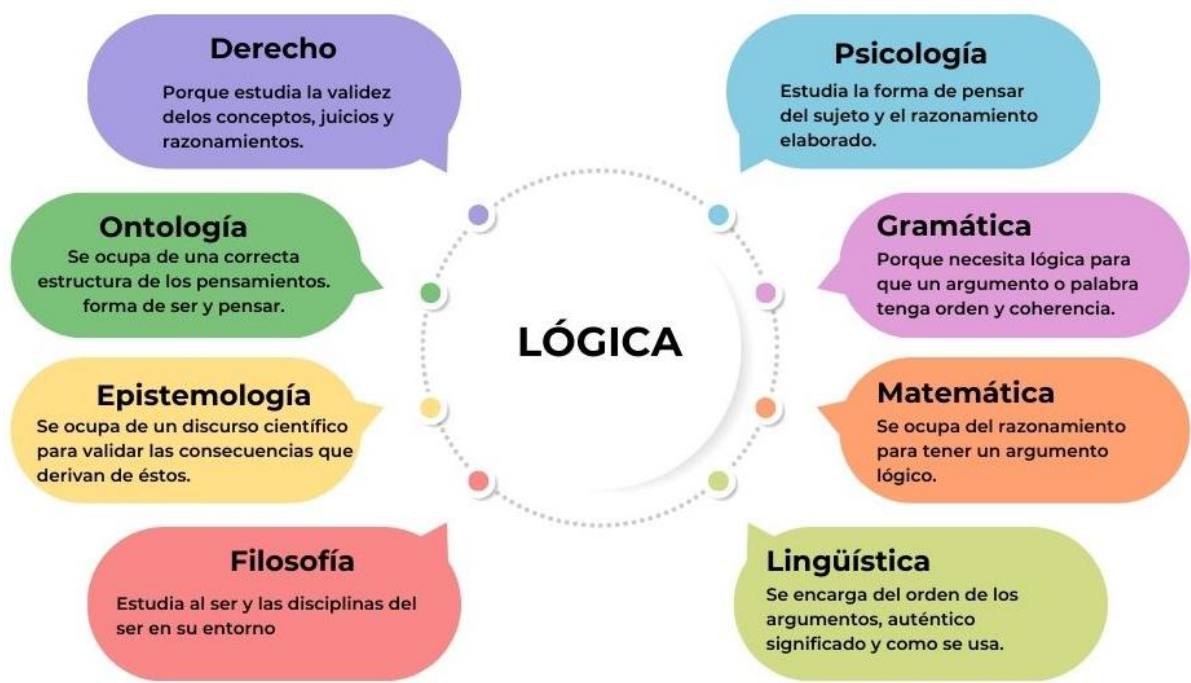

Nota. Adaptado de Tareas de Lógica y Apreciación Artística. por A. Aragón, 2013, Blog (http://themaschingon666.blogspot.com/p/unidad-i-logica.html). Todos los derechos reservados 2013 por Licenciatario.

## 2.2. Métodos lógicos para la demostración de argumentos.

Los métodos de demostración desde el punto de vista de la lógica consisten en un conjunto de procedimientos para llegar a la verdad o falsedad de los pensamientos. Algunos de ellos se describen brevemente a continuación, de acuerdo a autores expertos en el tema.

*2.2.1 Silogismo*

En palabras de Mateos (2016) es el raciocinio deductivo en el cual las premisas enlazan dos términos con un tercero, y la conclusión expresa la relación de estos dos términos entre sí.

*2.2.2 Tablas de verdad*

Como señalan Lozano y Pérez (2016) las tablas de verdad son un procedimiento que en una cantidad finita de pasos permite decidir mecánicamente la cantidad de semántica de una fórmula o fórmulas al revisar todas sus asignaciones relevantes de valor veritativo.

Este método de demostración gráfico que consta de filas y columnas permite interpretar la conclusión de un razonamiento y, aunque la clasificación según los valores obtenidos puede variar dependiendo de cada autor los más comunes son:

a) Tautología. Todos los valores son verdaderos.
b) Contingencia. Hay valores verdaderos y falsos.
c) Contradicción. Todos los valores son falsos.

Uno de los inconvenientes de este método es debido a que crecerá exponencialmente de acuerdo con el número de variables (enunciados) y por es considerado un algoritmo poco eficiente.

*2.2.3 Método directo*

Teniendo en cuenta a Tapia (2012) es el método donde se tienen como hipótesis verdaderas las proposiciones $P_1$, $P_2$, . . ., $P_n$ procediendo a la deducción de que la conclusión Q es verdadera a través de un proceso lógico deductivo.

*2.2.4 Método indirecto*

Tal como Tapia (2012) define este método como fundamento de equivalencia lógica entre las proposiciones $P \rightarrow Q$ y $\sim Q \rightarrow \sim P$. Para realizar una demostración contrapositiva se toma como hipótesis la negación de la conclusión escrita como $\sim Q$ para obtener como conclusión la negación de la hipótesis escrita como $\sim P$.

*2.2.5 Reglas de inferencia*

Hill (2013) refiere este método como un juego donde una vez aprendiéndose las reglas y su uso la validez consiste en el principio general de si las premisas son verdaderas, entonces las conclusiones derivadas de ellas lógicamente van a ser verdaderas.

## 2.3. Material didáctico mediado por las TIC

Como plantea Bates (1999), las TIC, traen consigo nuevas oportunidades para los entornos educativos, ya que posibilitan la inclusión de la diversidad a través de la combinación de medios que, aunados a un buen diseño instruccional, pueden ofrecer nuevas alternativas pedagógicas.

Gracias a las TIC los materiales didácticos mediados son una alternativa didáctica que cada día se encuentran de manera accesible también debido a la disponibilidad de los estudiantes con dispositivos móviles conectados a la red. En este sentido aun cuando en las TIC no está la solución de las dificultades del proceso de enseñanza-aprendizaje de las matemáticas si pueden provocar un cambio en la manera en que se enseñan, ya que mediante ellas se pueden plantear múltiples problemáticas permitiéndoles a los estudiantes desarrollar estrategias de solución y mejorar la comprensión de los conceptos matemáticos. A su vez permiten que los docentes se vuelvan facilitadores y busquen ofrecer desafíos y alternativas de trabajo a sus estudiantes, con el objetivo de ayudarles a construir su propio conocimiento.

Aunque existe una diversidad muy grande de tecnologías educativas a continuación se describirán las empleadas en el presente trabajo para poder elaborar el material didáctico y con ello lograr los objetivos establecidos en el temario de la asignatura de Estructuras Discretas:

Classroom. Permite crear una clase de manera fácil, distribuir y calificar tareas sin utilizar papel y comunicar con los alumnos en tiempo real.

Código QR. Quick Response Barcode. Esta tecnología permite cifrar, de forma rápida, texto plano en formato de código de barras. Los códigos QR conectan los objetos reales con cualquier contenido web.

Google Drive. Herramienta gratuita basada en Web para crear, subir y compartir documentos online con la posibilidad de colaborar en grupo.

Formularios de Google. Herramienta útil que permite planificar eventos, enviar una encuesta, hacer preguntas o recopilar otro tipo de información de forma fácil y sencilla.

Jamboard. Es una pizarra virtual en la que se puede escribir, dibujar, agregar imágenes y rayar sobre ellas.

Zoom. Plataforma que permite realizar videoconferencias, chatear e impartir clases de forma rápida y sencilla.

## 2.4. Modelo instruccional ASSURE

Fue ideado para ser utilizado por maestros en un salón de clases, para que los docentes o facilitadores puedan diseñar y desarrollar un ambiente de aprendizaje más apropiado para sus estudiantes. Este proceso puede utilizarse para planear las lecciones a enseñar y para mejorar la enseñanza y el aprendizaje.

Se basa en la teoría de aprendizaje cognitivista, no contempla una evaluación continua sino hasta el final, además de no realizar un análisis del entorno.

Según Heinich, et al., (1999) autores de este modelo, el acrónimo ASSURE consta de 6 fases que sirven de guía para poder planear y poder conducir la instrucción en el proceso de la enseñanza a distancia incorporando la tecnología de manera eficaz.

A continuación, se describirán las 6 etapas o fases De acuerdo con Benítez (2010):

La primera etapa consiste en *Analizar las características del estudiante* o de los participantes del curso, aspectos socioeconómicos y culturales, antecedentes escolares, edad, sexo, estilos de aprendizaje, etc. Estos aspectos ayudarán a realizar una mejor planeación.

La segunda etapa consiste en el *establecimiento de objetivos de aprendizaje,* de acuerdo con las características obtenidas de los estudiantes el establecimiento de los objetivos se debe especificar que comportamientos se van a evaluar, de tal manera que se la evaluación asegure el aprendizaje del estudiante.

La *selección de estrategias, tecnologías, medios y materiales* es la tercera etapa donde el docente debe realizar la selección para ser congruente con los objetivos de aprendizaje para facilitar que el estudiante logre los aprendizajes.

La cuarta etapa es la *Utilización de los medios y materiales*, donde consiste no solo en la implementación sino también el tener considerados otros medios a los seleccionados por si no llegarán a funcionar.

La quinta etapa denominada *participación de los estudiantes,* consiste en lograr que el estudiante comprenda y analice la información para lograr cumplir los objetivos establecidos mediante la participación del estudiante.

La ultima etapa es la *evaluación y revisión de la implementación y resultados del aprendizaje,* donde consiste en evaluar si se lograron los objetivos de aprendizaje haciendo uso de TIC.

## 2.5. Aprendizaje significativo

Para Ausubel (1963, p. 58), el aprendizaje significativo es el mecanismo humano, por excelencia, para adquirir y almacenar la inmensa cantidad de ideas e informaciones representadas en cualquier campo de conocimiento. El aprendizaje significativo va darse cuando una nueva información se une con un concepto relevante que ya existía en la estructura cognitiva del individuo.

El aprendizaje es la manera más importante por la cual una persona puede entender un tema, es decir promueve el aprendizaje significativo en lugar del aprendizaje de memoria.

El Modelo de Ausubel distingue tres tipos de aprendizajes que son:

• Aprendizaje de representaciones.

• Aprendizaje de conceptos.

• Aprendizaje de proposiciones.

Analizar las características de los estudiantes y del contexto en un diseño instruccional según Benítez (2010), para nosotros como docentes esta es la fase fundamental para que podamos obtener el aprendizaje significativo, ya que a partir de la información general y específica de nuestros estudiantes podremos elaborar de manera sencilla e imparcial nuestra planeación, sin dejar de tomar en cuenta el contexto donde se llevará a cabo.

Partir del análisis anterior nos permitirá realizar materiales y elegir las TIC (en determinado caso) adecuadas que cubran las verdaderas necesidades de nuestros estudiantes, de no llevarlo a cabo o no de la manera correcta podríamos perder tiempo y recursos, porque tendríamos que volver a rediseñar nuestra planeación.

# Capítulo 3. Diseño metodológico

## 3.1 Tipo de investigación

Se considera una investigación del tipo preexperimental, *cuasiexperimental, con grupos no equivalentes sólo con postest*, ya que se pretende establecer a los grupos a los que se imparte la materia de Estructuras Discretas en la Facultad de Ingeniería de la UNAM y poder realizar una comparación de un grupo que no recibe el material didáctico y otro que sí lo recibirá.

## 3.2 Sujetos de estudio

Se tuvo la participación de dos grupos de estudiantes de cuarto semestre de la carrera de Ingeniería en Computación de la Facultad de Ingeniería de la Universidad Nacional Autónoma de México cursando la materia de Estructuras Discretas. El rango de edades está entre los 20 y 21 años.

## 3.3 Técnicas e instrumentos

- Grupos no equivalentes sólo con postest.

Empleando las palabras de McMillan y Schumacher (2005), es un diseño donde se caracteriza por tener "un grupo de control o de comparación que no recibe tratamiento alguno o un tratamiento diferente se añade al diseño de un grupo solo con postest. En este proyecto el grupo de control fue el del semestre en curso 2021-1 y el segundo grupo fue el del semestre 2021-2.

- Cuestionario diagnóstico.

Como expresan Gardey y Pérez (2021), "un cuestionario es un conjunto de preguntas que se confecciona para obtener información con algún objetivo en concreto". En este caso en particular está constituido por 7 preguntas (ver Anexo 1), cuya finalidad fue aplicarlo previamente a los estudiantes y saber con qué herramientas tecnológicas contaban para poder tomar las clases.

- Cuestionario de opinión

Gerber (s.f) considera que un cuestionario de opinión es "un método de recolección de información, que, por medio de un cuestionario, recoge las actitudes, opiniones u otros datos de una población, tratando diversos temas de interés". En este estudio se emplearon dos, uno de ellos está constituido por 8 preguntas (ver Anexo 2) cuyo objetivo fue ver si el material y ejercicios resueltos durante las sesiones síncronas fueron los adecuados para lograr una buena comprensión y aprendizaje del tema de lógica proposicional y, el segundo consta de 15 preguntas (ver Anexo 3) con la finalidad de recoger información sobre si les fue de utilidad el material didáctico mediado con tecnología y sobre todo si les facilitó comprender mejor el tema de lógica proposicional.

## 3.4 Modelo de diseño instruccional

La intervención se basó en un diseño instruccional bajo el modelo creado por Heinich et al. (1993) llamado ASSURE dado que consta de seis etapas que permitieron realizar el análisis de los participantes, el establecimiento de los objetivos, la selección de los métodos, el uso de medios, así como de materiales, los requerimientos de participación de los estudiantes y por último como evaluar y revisar de acuerdo con el temario de la asignatura. Enseguida se detalla en qué consiste cada etapa.

*3.4.1 Etapa 1. Análisis de los participantes*

Participaron en el estudio estudiantes de la Facultad de Ingeniería, de la carrera de Computación, de la materia de Matemáticas Discretas, con edades entre 19 y 21 años, inscritos en dos grupos.

El primer grupo muestra fue el del semestre en curso 2021-1, el cual inicio el 21 septiembre 2020 y terminó el 29 enero 2020, con una población de 30 estudiantes -24 hombres y 6 mujeres-. El segundo grupo muestra fue el del semestre 2021-2, con una población de 38 estudiantes -29 hombres y 9 mujeres- de acuerdo con el calendario semestral de la UNAM inició el 22 febrero 2021 y terminaría el 19 de junio 2021, sin embargo, debido a que hubo un paro de actividades se recalendarizó el semestre y terminó hasta el 14 de agosto contemplando un período vacacional del 5 de julio 2021 al 24 de julio 2021.

*3.4.2 Etapa 2. Establecimiento de los objetivos*

La materia de matemáticas discretas tiene por objetivo que el alumno inferirá los conceptos matemáticos mediante la computación en la solución de problemas relacionados con el procesamiento de la información, el diseño de computadoras y de programas.

Objetivo de la lógica proposicional es que el alumno infiera la teoría de la lógica matemática para aplicarla en la solución de problemas dentro del campo de la computación

*3.4.3 Etapa 3. Selección de medios y materiales*

En la Tabla 1 se muestra de acuerdo con el número de sesiones los medios y materiales empleados para cada semestre, teniendo en cuenta que para el semestre 2021-1 el objetivo era no aplicar el material didáctico mediado por TIC.

**Tabla 1**

*Medios y materiales seleccionados para los semestres 2021-1 y 2021-2*

| | | **Medios y materiales** | |
|---|---|---|---|
| | **No. sesiones** | **Medios** | **Materiales** |
| **Semestre 2021-1** | 7 | Pizarra virtual de Zoom | Documentos digitales en formato pdf |
| | 1 | Aplicación Zoom<br>Correo electrónico | Documentos digitales en formato pdf |
| **Semestre 2021-2** | 5 | Aplicación Zoom, Jamboard, Google Classroom, Google Drive, Código QR | Documentos digitales en formato pdf |
| | 1 | Aplicación Zoom<br>Correo electrónico | Documentos digitales en formato pdf |

### *3.4.4 Etapa 4. Utilización de medios y materiales*

Como se puede observar en la Tabla 2 la intervención para el semestre 2021-1 (del 21 septiembre 2020 al 29 enero 2021) se llevó a cabo en ocho sesiones, donde el objetivo de las primeras siete, fue emplear recursos tecnológicos básicos que ayudan a recrear lo más parecido posible a como se venían realizando las clases de forma presencial y poder justificar la creación del material didáctico con el uso de tecnología educativa *y con ello* contribu*ir* en el logro de los objetivos establecidos en el temario de la asignatura.

Cada sesión fue de aproximadamente 110 minutos de manera síncrona, los lunes y miércoles mediante la herramienta de videoconferencias Zoom. Se dio una explicación teórica y práctica sobre el tema de lógica proposicional empleando la pizarra que trae integrada Zoom. Al término de la sesión se compartió dicha pizarra a los estudiantes en un documento pdf.

En la octava sesión se aplicó un examen con una duración de noventa minutos que se les hizo llegar vía correo electrónico, cuyas respuestas debían ser enviadas por el mismo medio.

**Tabla 2.**

| | Fase 2. Diseño de la intervención |
|---|---|
| Herramientas utilizadas. | Jamboard, Zoom, Google Forms, Telegram, Gmail, dispositivo móvil o Tablet (mlearning), Google Drive y Google Classroom. |
| Perfil de los participantes. | Facilitador: debe conocer la disciplina. Capacidad para coordinar, orientar, potenciar la autonomía, el trabajo cooperativo y la toma de decisiones en el estudiante. Flexibilidad en el seguimiento del proceso formativo y propiciar la interacción entre los estudiantes.<br>Estudiante: debe de contar con habilidades cognitivas de abstracción, de análisis, de síntesis y de reflexión. Capacidad de comunicación oral y escrita empleando correctamente el lenguaje de las matemáticas. |
| Contenidos. | • Fórmulas proposicionales y tablas de verdad.<br>• Formas normales y dispositivos de dos estados.<br>• Notación Polaca y parentizada.<br>• Elementos de inferencia para el cálculo proposicional.<br>• Prueba automática de teoremas.<br>• Fórmulas de predicados. |
| Estrategias y técnicas de enseñanza. | Explicación teórica y práctica por parte del facilitador.<br>Uso de la plataforma de Google Classroom.<br>Ejercicios con código QR que dirigen a un video alojado en Google Drive. |
| Recursos (videos, podcast, texto, entre otros). | Ejercicios en formato digital y videos. |
| Comunicación (facilitador-participante, participante-participante, entre otros). | Comunicación síncrona: una retroalimentación del tema (facilitador-participante), trabajo en equipo por parte de los estudiantes (participante-participante). Comunicación asíncrona mediante correo electrónico, Google Classroom y Telegram se podrán plantear dudas por parte del participante y resueltas por el facilitador. |

A continuación, en las tablas 3 y 4 se detalla la estrategia de intervención de los semestres 2021-1 y 2021-2 respectivamente.

**Tabla 3**

*Estrategia de intervención del semestre 2021-1*

| No. sesiones | Tiempo de total (min) | Fecha de aplicación | Actividad | Recursos | Procedimiento | Objetivo |
|---|---|---|---|---|---|---|
| 7 | 770 | 28/09/2020 30/09/2020 05/10/2020 07/10/2020 12/10/2020 14/10/2020 19/10/2020 | Clase síncrona a distancia | Pizarra virtual de Zoom Documentos digitales en formato pdf | Explicación teórica y práctica sobre el tema de lógica proposicional. | Impartir clases de manera síncrona a distancia empleando las herramientas tecnológicas mínimas requeridas. |
| 1 | 90 | 21/10/2020 | Examen | Documentos digitales en formato pdf Correo electrónico | Envío y recepción de un documento digital en formato pdf vía correo electrónico | Evaluar el nivel de comprensión del estudiante sobre el tema de lógica proposicional |

En la Tabla 3 se expone la intervención para la segunda muestra que fue la del semestre 2021-2 (del 22 febrero 2021 al del 19 de junio 2021), esta se llevó a cabo en seis sesiones donde se logró impartir las primeras cinco clases de aproximadamente 110 minutos cada una de manera síncrona mediante la aplicación Zoom además de complementar el uso de la plataforma de Google Classroom donde se subieron documentos digitales con un código QR que direccionaban a un video alojado en Google Drive para que el estudiante los pudiera consultar en el momento que lo deseara. Sus dudas, comentarios u observaciones también se realizaron por el mismo medio.

El objetivo de dichas sesiones fue emplear recursos tecnológicos educativos para llevar a cabo las clases de manera síncrona con la implementación del material didáctico mediado por tecnología para poder facilitarles el aprendizaje sin que interfiriera en el tiempo de las clases síncronas.

Para la sexta sesión de aplicó un examen de noventa minutos con el fin de evaluar el nivel de comprensión del estudiante sobre el tema de lógica proposicional. Nuevamente se les hizo llegar por vía correo electrónico y posteriormente los estudiantes enviaron las respuestas por el mismo medio.

**Tabla 4**

*Estrategia de intervención del semestre 2021-2*

| No. sesiones | Tiempo de total (min) | Fecha de aplicación | Actividad | Recursos | Procedimiento | Objetivo |
|---|---|---|---|---|---|---|
| 5 | 550 | 26/07/2021 28/07/2021 02/08/2021 04/08/2021 09/08/2021 | Clase síncrona a distancia Implementación de material didáctico mediado por tecnología | Aplicación Zoom Jamboard Documentos digitales en formato pdf Google Classroom Google Drive Lector código QR | Explicación teórica y práctica sobre el tema de lógica proposicional. En la plataforma Google Classroom se subieron documentos digitales con un código QR que direccionarán a un video alojado en Google Drive | Implementar material didáctico mediado por tecnología para facilitar el aprendizaje a los estudiantes sobre el tema de lógica proposicional sin que interfiera en el tiempo de las clases síncronas. |
| 1 | 90 | 21/10/2020 | Examen | Documentos digitales en formato pdf Correo electrónico | Envio y recepción de un documento digital en formato pdf vía correo electrónico | Evaluar el nivel de comprensión del estudiante sobre el tema de lógica proposicional |

*3.4.5 Etapa 5. Participación de los estudiantes*

Para los estudiantes del semestre 2021-1 sólo se diseñaron ejercicios, tareas, exámenes, material didáctico y con el apoyo de herramientas tecnológicas que hicieron crear un escenario lo más cercano a lo presencial.

Para los estudiantes del semestre 2021-2 se aplicó básicamente el mismo procedimiento, pero esta vez se incorporó el uso de herramientas que pudieron emplear no solo desde su computadora personal ahora también tuvieron la oportunidad de visualizarlo desde su dispositivo móvil o Tablet (mlearning) para poder tener acceso a los documentos digitales dentro y fuera de clase, con ejemplos resueltos paso a paso y, problemas propuestos en la plataforma de Google Classroom, desde esta misma plataforma pudieron exponer sus dudas o comentarios manteniendo una comunicación asíncrona. Además, los ejemplos resueltos contaron con un código QR que los dirigió a un video alojado en Google Drive donde contiene la explicación detallada de los ejercicios.

*3.4.6. Etapa 6. Evaluación*

La evaluación se llevó a cabo para ambos semestres mediante ejercicios, tareas y un examen de conocimiento.

Una vez obtenidas las calificaciones y promedios de ejercicios, tareas y exámenes, se graficaron para compararlos los resultados del semestre 2021-1 con los resultados del semestre 2021-2.

También se aplicó un cuestionario de opinión para determinar si los estudiantes consideraban necesario más ejercicios para reforzar lo aprendido y con ello haber obtenido un mejor resultado, esto con el objetivo de ver desde el punto de vista de los estudiantes si realmente el material y técnicas que se emplearán en el semestre 2021-2 eran requeridas o no. El resultado de este cuestionario también se graficó.

Se aplicó a los estudiantes del semestre 2021-2 un cuestionario de opinión para obtener su juicio acerca de si les fue útil y apropiado el material proporcionado para mejorar su aprendizaje de lógica proposicional de acuerdo con los objetivos de la asignatura. Estos resultados también se graficaron.

# Capítulo 4. Resultados

Inicialmente de acuerdo con el diseño metodológico la muestra del semestre 2021-1 (del 21 septiembre 2020 al 29 enero 2021) debido a que el primer y el único contacto con los estudiantes era a través de su dirección electrónica, por esto se empleó la herramienta Gmail de Google como servicio de correo electrónico para crear un alias y poder establecer comunicación con ellos. Dadas dichas circunstancias se diseñó y planeó un cuestionario diagnóstico elaborado con Google Forms (ver Anexo 1), para conocer los recursos tecnológicos con que contaban y poder planear la mecánica de las sesiones así como la comunicación durante el semestre.

Como se puede observar en la Figura 2, al momento de levantar la información, los 30 estudiantes disponían con una computadora o laptop y 28 de ellos con un dispositivo móvil. A pesar de que sólo 19 de ellos tenían escáner, podían instalar una aplicación en su móvil para poder enviar sus tareas y ejercicios.

**Figura 2**

*Dispositivos con los que cuenta el estudiante en casa.*

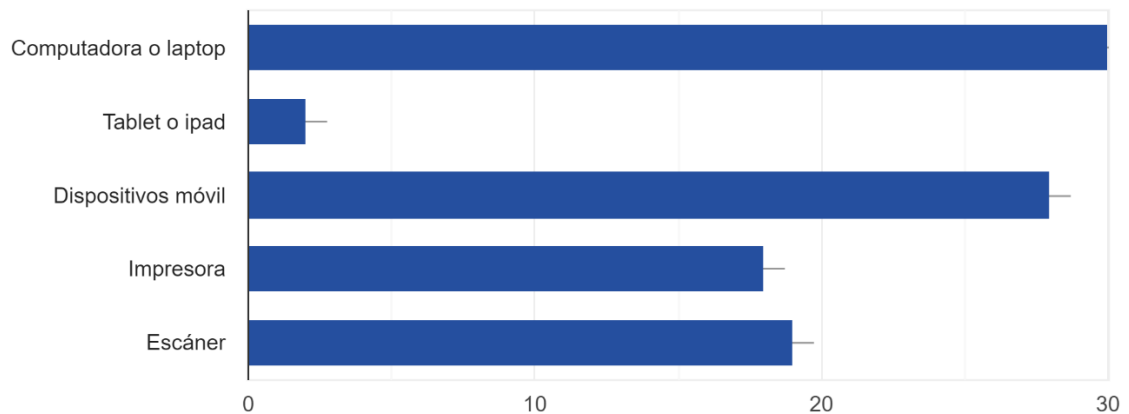

En la Figura 3, se puede observar que de los 30 estudiantes sólo uno empleaba internet por medio de un plan de datos, el resto lo hacía mediante Wi-Fi, por ello se agregó la herramienta de videoconferencias Zoom para las clases síncronas, además la UNAM proporcionó a todos los docentes una licencia para poder emplearla sin tener restricción de tiempo.

**Figura 3**

*Disponibilidad de internet del estudiante en casa*

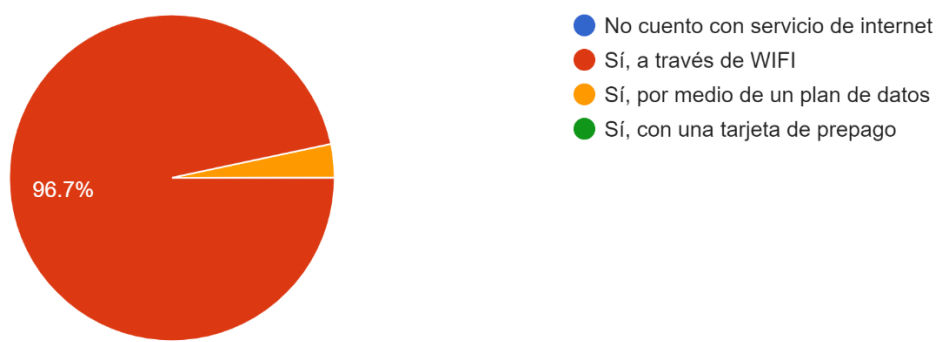

**Nota:** Los estudiantes tenían cuatro opciones para responder si contaban o no con servicio de internet y a través de que medio lo obtenían.

En la Figura 4, se muestra que no todos los estudiantes podían estar conectados, por lo que se añadió la herramienta de Telegram para mantener una comunicación asíncrona sin invadir su privacidad al tener que solicitarles su número telefónico y respetar el aviso de privacidad establecido no sólo en la Facultad de Ingeniería sino también a nivel UNAM.

**Figura 4**

*Tiempo de conexión de los estudiantes de forma síncrona*

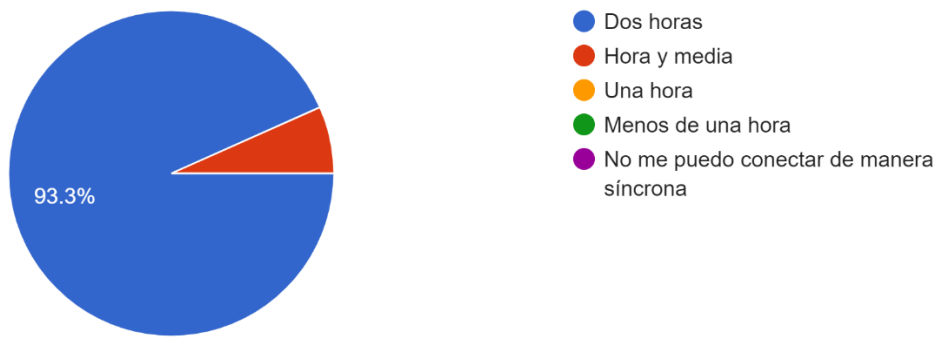

Con relación a las cuentas de correo electrónicos los estudiantes contaban con diferentes dominios y el trámite de correos institucionales es tardado, asimismo coincidió con una eventualidad que tuvo el área de TI de la Facultad de Ingeniería con Google, por lo que no fue posible emplear Google Classroom. Sin embargo, se utilizó como medio por el que se envió y recibió tareas, ejercicios y lo visto en clase.

Inicialmente esta intervención había sido planeada de manera presencial y la aplicación del uso de material didáctico mediado con herramientas tecnológicas sería hasta el semestre 2021-2, sin embargo, por las circunstancias de la pandemia se tuvo que realizar material didáctico y emplear algunas herramientas para poder dar la clase a distancia, pero tratando de crear un escenario lo más cercano a lo presencial.

Previo a cada clase se diseñaron y elaboraron ejercicios, tareas y material digital complementario empleando Microsoft Word (ver Anexo 4) para poder enviarlo por correo electrónico en formato pdf y evitar problemas en cuanto a las diferentes versiones del programa o sistema operativo que cada estudiante empleaba.

Se puso en práctica la técnica del aula invertida por lo que se les proporcionó la teoría en formato digital por correo electrónico para ganar tiempo y evitar estar conectados por videoconferencia demasiado tiempo.

Durante las clases síncronas se empleó la pizarra de Zoom que es como usar un pizarrón de manera presencial, además de que permite guardarlas y posteriormente compartirlas con los estudiantes en formato pdf.

Se pudo intervenir como se tenía planeado, aunque hubo muchos problemas de conexión a Internet aunado a que los estudiantes tuvieron dudas de temas anteriores que no eran de lógica proposicional las cuales fueron aclaradas durante las sesiones síncronas de Zoom.

Para la segunda muestra correspondiente al semestre 2021-2 con fecha de inicio del 22 febrero 2021 y como fecha final del 19 de junio 2021, no se pudo llevar a cabo debido al paro de actividades promovido por las estudiantes en apoyo a los profesores de asignatura que no estaban recibiendo sus honorarios correspondientes. Dicha situación desencadenó varias modificaciones como fue la recalendarización del semestre, reducción del número de sesiones establecidas originalmente en el temario de la asignatura, libertad por parte del docente para establecer el número de sesiones para cada tema, así como los temas a abordar. Lo único que se debía cumplir era terminar el 14 de agosto del 2021 el semestre, contemplando un período vacacional del 5 de julio 2021 al 24 de julio 2021. Con respecto a la que la primera muestra (semestre 2021-1), hubo una diferencia de dos sesiones menos de las consideradas en estrategia de intervención.

Tanto el primero como el segundo objetivo específicos se cumplieron al cien por ciento, ya que la selección de las herramientas de software adecuadas para implementación y producción del material didáctico no solo estuvieron enfocadas para poder cumplir con los objetivos de la asignatura, sino también a las herramientas con licencia e institucionales que la UNAM puso a disposición de los docentes – herramientas y recursos de *Google*, *Zoom*, *Webex*, *Moodle* -. En el Anexo 5 se presenta el rediseño de 15 ejercicios de lógica proposicional por diferentes métodos de solución de forma detallada y desglosando paso a paso su solución, empleando para ello *Jamboard* y Documentos de *Google*, guardados en formato pdf, así como en el mismo *Jamboard*.

La grabación de 14 videos se realizó con *Zoom* alojados en *Google Drive* y con el generador gratuito *qrcodemokey* los códigos QR para poder acceder a ellos.

En la Tabla 4 se observa que del tema de Formas Normales se rediseñaron 5 ejercicios cada uno por un solo método de solución. Del tema Directo e Indirecto fueron 4 ejercicios con un total de 5 soluciones distintas y cada ejercicio por 2 métodos diferentes. Del tema de Reglas de Inferencia, se rediseñaron 6 ejercicios donde fueron un total de 33 soluciones distintas y cada uno de ellos se resolvió por 2 métodos diferentes.

**Tabla 4**

*Relación de ejercicios con su número de soluciones*

| Tema | Número de ejercicios | Métodos de solución diferentes | Número de soluciones diferentes |
|---|---|---|---|
| Formas Normales | 5 | 1 | 5 |
| Directo e indirecto | 4 | 2 | 5 |
| Reglas de inferencia | 6 | 2 | 33 |

Para poder llevar a cabo el tercer objetivo específico sobre la implementación de los materiales didácticos en apego a los objetivos de la asignatura, se realizó el taller "Enseñar y evaluar con enfoque formativo" impartido en la Coordinación de Universidad Abierta, Innovación Educativa y Educación a Distancia de la UNAM, con la finalidad de crear un aula digital de Google Classroom y mediante ella poner a disposición de los estudiantes el material didáctico. Dicha aula fue revisada, evaluada y aprobada por el Mtro. Alejandro González Flores instructor del taller (ver Anexo 6).

En cuanto al cuarto objetivo específico referente al análisis del rendimiento académico después de haber implementado el material didáctico, se graficaron para ambas muestras las calificaciones de un examen, así como las calificaciones promedio de ejercicios realizados de manera independiente por los estudiantes.

En la Figura 5 se presenta la gráfica de las calificaciones del examen correspondiente al semestre 2021-1, donde 16 de los 29 estudiantes obtuvieron una calificación igual o mayor a 9, 11 estudiantes con una calificación entre 7 a 8.9, un solo estudiante con calificación entre 6 a 6.9 y, solo uno con calificación reprobatoria. Obteniendo un promedio general de 8.8.

**Figura 5**

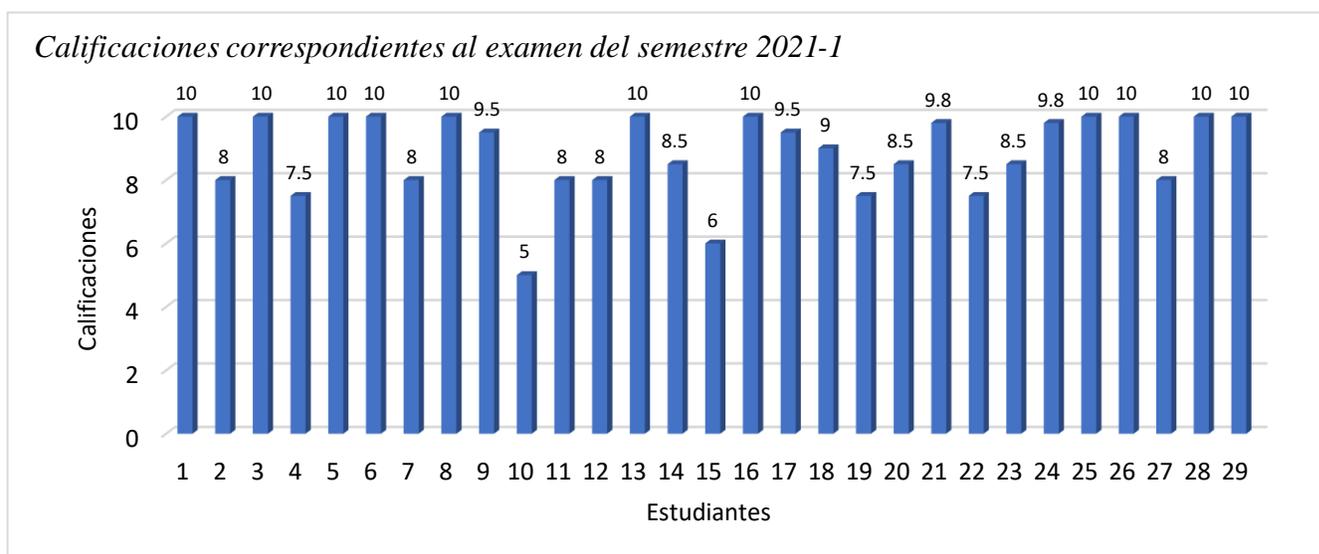

*Calificaciones correspondientes al examen del semestre 2021-1*

Y en la Figura 6, se pueden observar las calificaciones del examen correspondiente al semestre 2021-2, de los 38 estudiantes 19 obtuvieron una calificación igual o mayor a 9, 8 estudiantes con una calificación de 7 a 8.9, 5 estudiantes con calificación de 6 a 6.9 y, seis con calificación reprobatoria. Obteniendo un promedio general de 8.0.

**Figura 6**

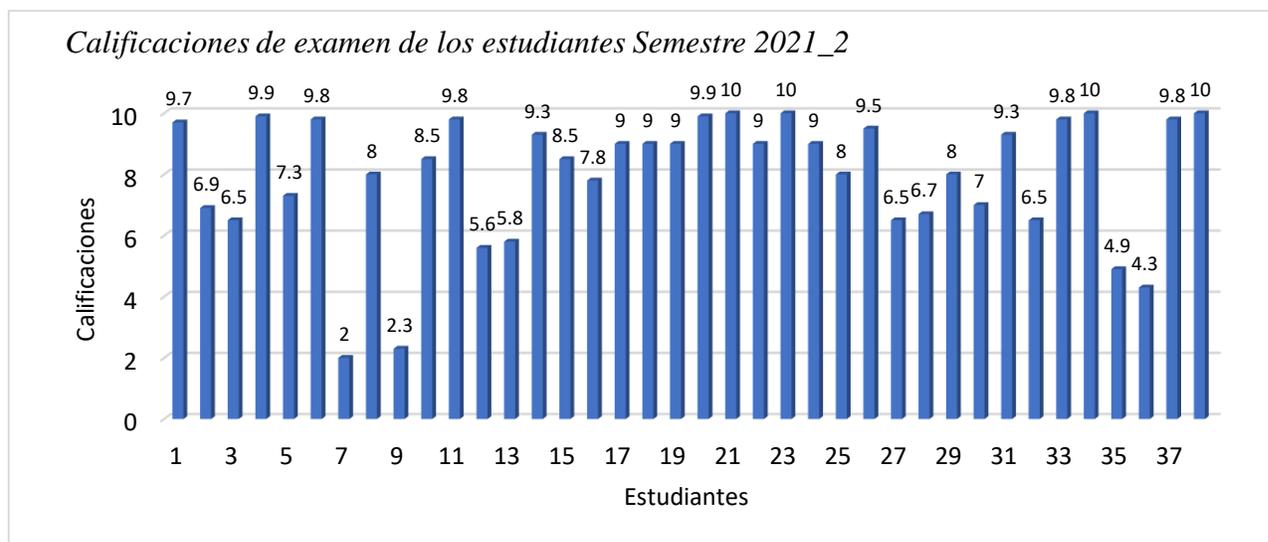

*Calificaciones de examen de los estudiantes Semestre 2021_2*

De acuerdo con el número de estudiantes de cada muestra se observa que el 50% su calificación fue mayor o igual a 9, a pesar de haber una diferencia mayor de 0.8 en cuanto al promedio general cuando los estudiantes no emplearon el material didáctico se debe considerar que el número de sesiones síncronas fue menor debido al paro de labores de la Facultad.

Por otro lado, en la Figura 7, se muestran las calificaciones promedio correspondientes al semestre 2021-1. Los estudiantes realizaron de manera independiente seis ejercicios, donde 7 (que representan el 24.1%) de los 29 estudiantes obtuvieron una calificación igual o mayor a 9, con una calificación promedio general de 8.3.

**Figura 7**

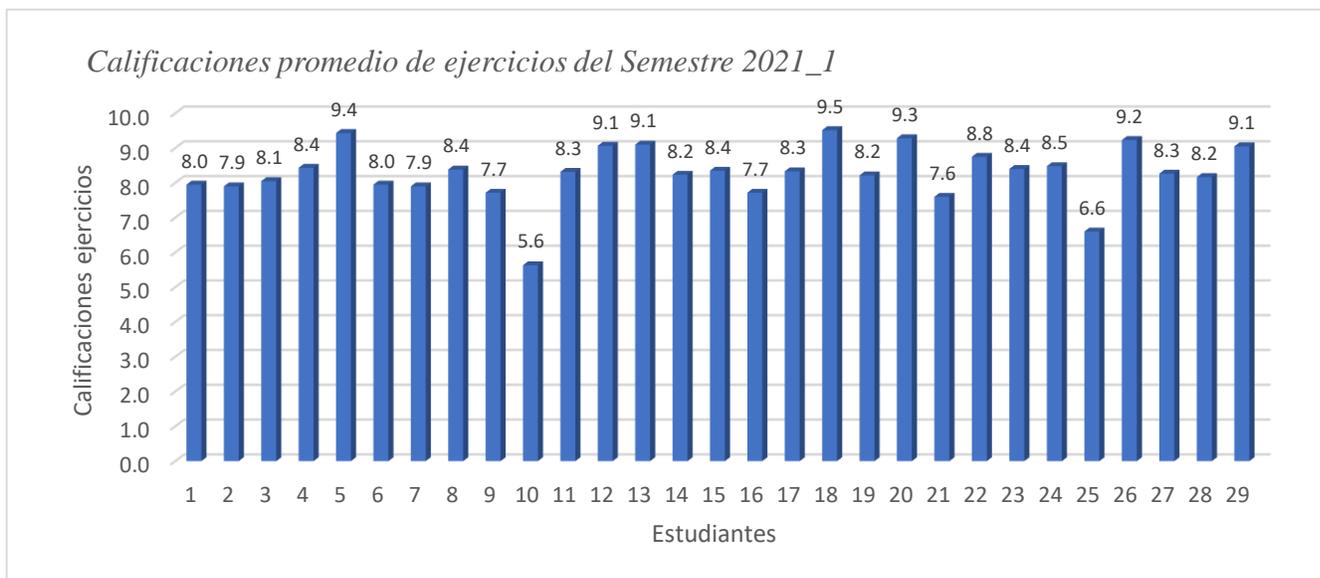

Para el semestre 2021-2 solo fue posible que realizaran tres ejercicios de manera independiente debido a la restructuración del semestre. Las calificaciones promedio se pueden observar en la Figura 8, de los 38 estudiantes 17 (representan el 44.7%) obtuvieron una calificación igual o mayor a 9 y un promedio general de 8.5.

**Figura 8**

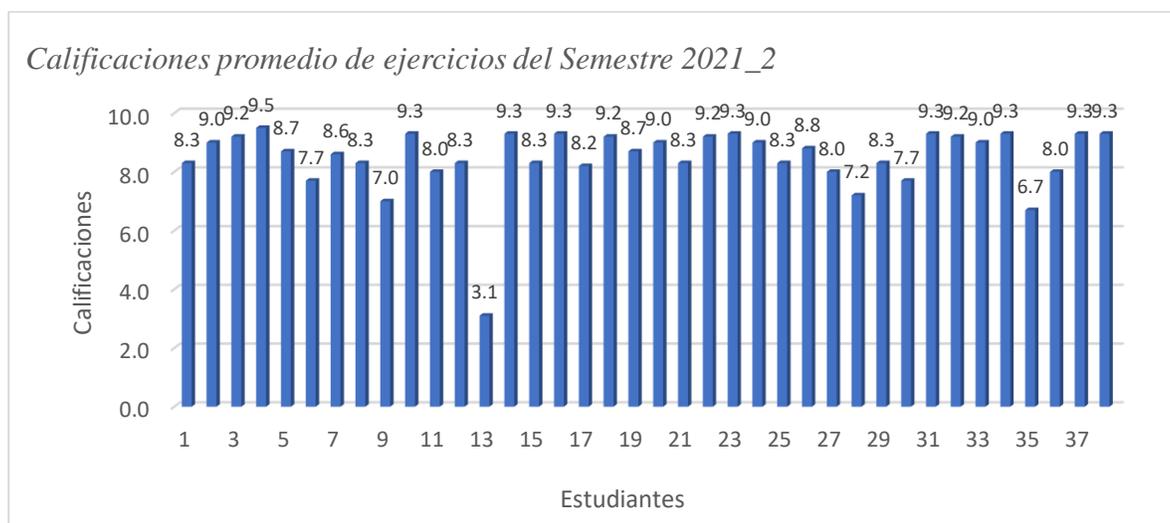

De tal manera que en la segunda muestra si hay una diferencia positiva del 20.6 % en cuanto a las calificaciones mayores o iguales a 9 y de 0.2 más con respecto al promedio general.

Para finalizar con el cumplimiento de los objetivos específicos en cuanto a medir el nivel de satisfacción del curso mediado por tecnología, se le aplicó el cuestionario de opinión a la primer muestra (ver Anexo 2). Como se puede observar en la Figura 9, el 60.7 % de los estudiantes consideran que son necesarios más ejercicios resueltos paso a paso como material de apoyo para una mejor comprensión del tema de lógica proposicional.

**Figura 9**

*Necesidad de material de apoyo con ejercicios resueltos*

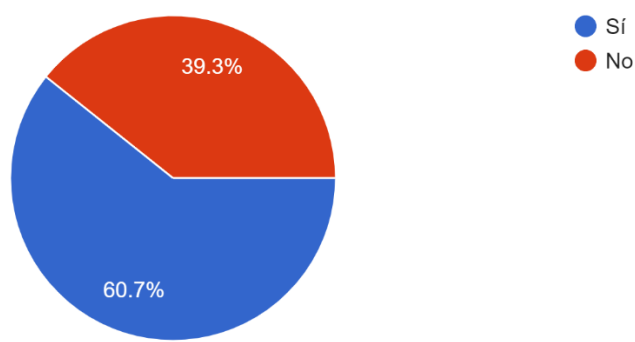

Nota: Esta figura muestra la opinión de los estudiantes del semestre 2021-1 en cuanto a la necesidad de la creación de material de apoyo con ejercicios resueltos paso a paso.

De manera similar en la Figura 10, el 64.3% de los estudiantes consideran que la explicación de ejercicios mediante videos ayudaría a mejorar la comprensión del tema de lógica proposicional.

**Figura 10**

*Ejercicios explicados en videos ayudarían en la compresión de lógica proposicional*

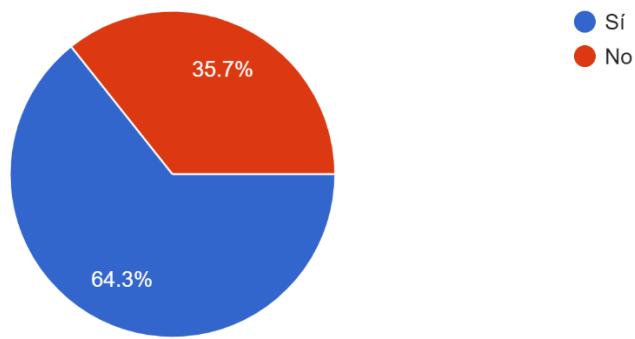

Posterior a la realización del examen del semestre 2021-2 (segunda muestra), se aplicó la encuesta de opinión (ver Anexo 3), en donde desde el punto de vista de los estudiantes el 73.7% de ellos consideran insuficientes el número de ejercicios resueltos en clase, igualmente, el 84.2 % consideró que el material complementario con ejercicios resueltos fue suficiente. Dichos resultados se exponen en las Figuras 11 y 12, respectivamente.

**Figura 11**

*Suficientes el número de ejercicios realizados por el docente durante la clase*

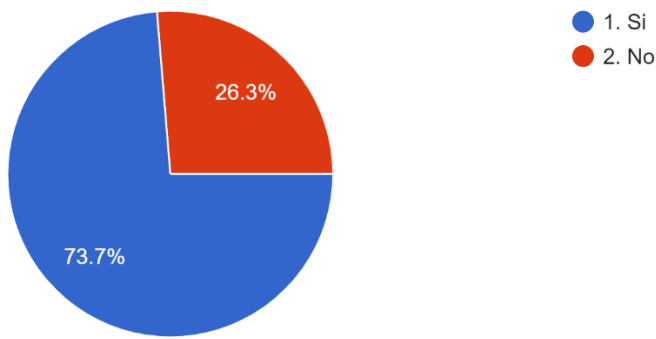

Nota: Esta figura muestra la opinión de los estudiantes del semestre 2021-2, donde el 73.7% consideró suficiente el número de ejercicios realizados por el docente durante la clase.

**Figura 12**

*Material complementario con ejercicios resueltos y explicados en video.*

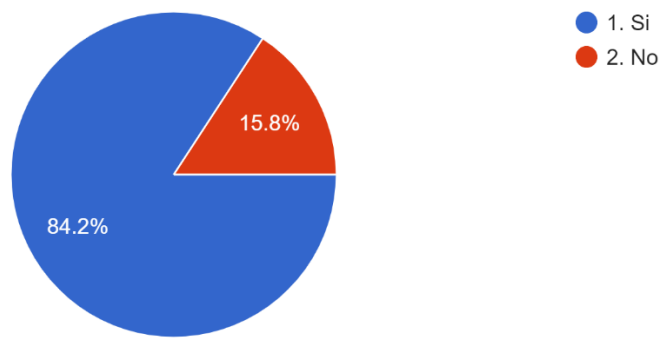

Nota: Esta figura muestra la opinión de los estudiantes del semestre 2021-2, donde el 84.2% consideró suficiente el material complementario con ejercicios resueltos y explicados en video.

# Capítulo 5. Conclusiones

Las técnicas elegidas junto con el diseño de la intervención guiaron la realización de una planificación adecuada en la obtención de los resultados y poder analizar si se logró cumplir el objetivo principal, así como los específicos enfocados principalmente en facilitar el aprendizaje de la lógica proposicional en los estudiantes de la asignatura de Estructuras Discretas de la Facultad de Ingeniería de la Universidad Nacional Autónoma de México.

En relación con el cumplimiento de los objetivos específicos como se pudo observar en los resultados obtenidos se cumplieron en su totalidad. Se pudo rediseñar, producir e implementar los materiales didácticos mediados por tecnología en apego a los objetivos de la asignatura. También se pudo analizar el rendimiento académico después de haber implementado el material didáctico; además cabe subrayar que el nivel de satisfacción conforme a los resultados de los cuestionarios de opinión aplicados a los estudiantes si reflejaron de manera cualitativa un porcentaje del más del 50% por una respuesta positiva, es decir, el primer cuestionario justificó la implementación de las herramientas y del material didáctico, asimismo el segundo reflejó la utilidad de la implementación en beneficio de una mejora en el aprendizaje de la asignatura.

Por otro lado, aunque el análisis cuantitativo mediante las gráficas de las calificaciones obtenidas no reflejó una mejora significativa en el semestre 2021-2 donde se implementó el material didáctico con el uso de tecnología, hubo dos factores no controlables e inesperados como la pandemia y el paro de actividades de la Facultad de Ingeniería, provocando dos escenarios completamente diferentes en cada una de las muestras de tal manera que, de no haber sido por la implementación del material probablemente se hubieran visto afectados significativamente los estudiantes en su aprendizaje de la lógica proposicional.

Por tanto, se puede decir que se logró el objetivo general de este estudio que consistió en facilitar el aprendizaje de la lógica proposicional mediante la implementación de material didáctico con el uso de tecnología educativa puesto que sí contribuyó en el logro de los objetivos establecidos en el temario de la asignatura.

En futuros estudios lo ideal es que las muestras sean dos grupos de Estructuras Discretas en el mismo semestre para que el análisis tanto cualitativo como cuantitativo de los resultados obtenidos estén bajo las mismas características y factores, de tal manera que se puedan realizar las adecuaciones necesarias, ya que todo es perfectible.

Considerando lo anterior es recomendable también tener prácticas prediseñadas, pero con la flexibilidad de poder crear similares con mayor o menor grado de dificultad de acuerdo con el desempeño y resultado de los estudiantes.

Por otro lado este proyecto es el comienzo de una tarea ardua para poder lograr un compendio de materiales que les sea de utilidad a los estudiantes de las carreras enfocadas al campo informático, ya que , aunque en un principio se pensó únicamente en un solo grupo de la Facultad de Ingeniería, el beneficio se puede extender a los demás grupos donde se imparte la materia y porque no pensar en otras instituciones.

# Referencias

# Anexos 57## Anexo 1

*Instrumento de diagnóstico*

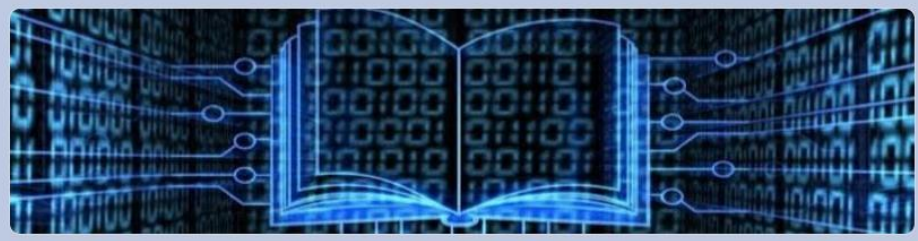

¿Cuántas personas ocupan en la misma casa, el servicio de Internet para sus actividades escolares o laborales? *

- ◯ Una
- ◯ Dos
- ◯ Tres
- ◯ Más de cuatro

En caso de tener computadora, laptop o tablet, selecciona la opción que mejor describa el tiempo que puedes ocupar, diariamente dicho dispositivo para tus actividades escolares. *

- ◯ Todo el tiempo necesario ya que es para mi uso personal y exclusivo.
- ◯ El tiempo suficiente, ya que la comparto con otra persona, pero nos organizamos.
- ◯ Un tiempo limitado, ya que la comparto con dos o más personas.
- ◯ Casi nada de tiempo, porque hay una sola computadora y la usamos muchas personas.

En el horario de clase ( 7:00 a 9:00 hrs) , ¿cuánto tiempo puedes estar conectado? *

- ◯ Dos horas
- ◯ Hora y media
- ◯ Una hora
- ◯ Menos de una hora
- ◯ No me puedo conectar de manera síncrona

¿Tienes un lugar específico en tu casa habilitado para hacer tus trabajos y/o tareas? *

- ◯ Sí
- ◯ No
- ◯ No, pero puedo realizar mis tareas sin interrupciones.

**Nota:** El instrumento de diagnóstico fue un cuestionario elaborado con Forms de Google. El cual fue aplicado a los estudiantes de ambos semestres.

De acuerdo con los resultados, se puede observar que es posible realizar clases de manera síncrona o hacer uso de herramientas tecnológicas que requieran el uso de internet para estar llevando a cabo sus actividades académicas.

**Anexo 2**

*Cuestionario de opinión 1*

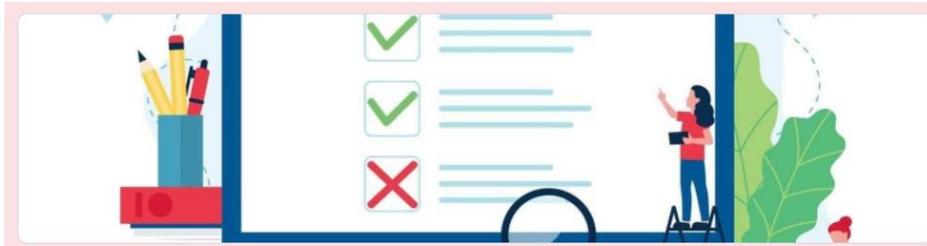

4. ¿El número de tareas fue el adecuado? *

○ Sí

○ No, fueron pocas

○ No, fueron muchas

5. ¿El material complementario con ejercicios resueltos proporcionado por la profesora fue suficiente? *

○ Sí

○ No

6. ¿El nivel de dificultad de los ejercicios resueltos y explicados por la profesora fue el adecuado para que tú pudieras resolver otros por tu cuenta? *

○ Si

○ No

7. ¿Consideras necesario más ejercicios resueltos paso a paso como material de apoyo para una mejor comprensión del tema? *

○ Sí

○ No

8. ¿Consideras que la explicación de ejercicios mediante videos ayudaría a mejorar la comprensión del tema? *

○ Sí

○ No

**Anexo 3** 62

*Cuestionario de opinión 2*

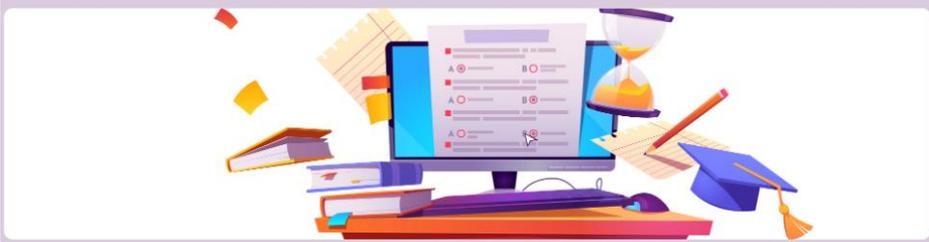

Semestre: *

- ◯ 1. Cuarto
- ◯ 2. Quinto
- ◯ 3. Sexto

Número de veces que has cursado la materia: *

- ◯ 1. Primera
- ◯ 2. Segunda
- ◯ 3. Tercera
- ◯ 4. Cuarta

Sección 3 de 4

## Material y ejercicios proporcionados.

En esta sección se solicita información sobre el material y ejercicios proporcionados.

¿Para comprender el tema basta con a explicación de la profesora? *

- ◯ 1. Si
- ◯ 2. No

¿Consultaste TODO el material proporcionado por la profesora? *

- ◯ 1. Si
- ◯ 2. No

¿Consideras que el material proporcionado con la teoría es adecuado? *

○ 1. Si

○ 2. No

¿El número de ejercicios realizados por la profesora en clase fueron suficientes?

○ 1. Si

○ 2. No

¿El número de ejercicios realizados por el estudiante en clase fueron suficientes? *

○ 1. Si

○ 2. No

¿El número de tareas fue el adecuado? *

○ 1. Si

○ 2. No

¿El material complementario con ejercicios resueltos proporcionado por la profesora fue suficiente? *

○ 1. Si

○ 2. No

¿El nivel de dificultad de los ejercicios resueltos y explicados por la profesora fue el adecuado para que tú pudieras resolver otros por tu cuenta? *

◯ 1. Si

◯ 2. No

¿Consideras necesario más ejercicios resueltos paso a paso como material de apoyo para una mejor comprensión del tema? *

◯ 1. Si

◯ 2. No

¿Consideras que la explicación de ejercicios mediante videos ayudaría a mejorar la comprensión del tema? *

◯ 1. Si

◯ 2. No

¿Consideras que la explicación de ejercicios mediante material impreso ayudaría más a mejorar la comprensión del tema que mediante videos? *

◯ 1. Si

◯ 2. No

Sección 4 de 4

## Material de consulta.

En esta sección se solicita información sobre el tipo de material que se consulta para complementar lo visto en clase.

Favor de seleccionar según corresponda tu respuesta: *

|  | Nunca (1) | Casi nunca (2) | Regular (3) | Casi siempre (4) | Siempre (5) |
|---|---|---|---|---|---|
| Libros impresos | ○ | ○ | ○ | ○ | ○ |
| Libros publicad… | ○ | ○ | ○ | ○ | ○ |
| Publicaciones c… | ○ | ○ | ○ | ○ | ○ |
| Publicaciones c… | ○ | ○ | ○ | ○ | ○ |
| Páginas web | ○ | ○ | ○ | ○ | ○ |
| Vídeos | ○ | ○ | ○ | ○ | ○ |
| Archivos digital… | ○ | ○ | ○ | ○ | ○ |
| Presentaciones | ○ | ○ | ○ | ○ | ○ |
| Infografías | ○ | ○ | ○ | ○ | ○ |
| Imágenes | ○ | ○ | ○ | ○ | ○ |

# Anexo 4

**FORMAS NORMALES**

**Literal, literales complementarios, par complementario.** Un literal es un átomo o la negación de un átomo. Si P, Q, R son átomos P, ¬P, Q, ¬Q, R, ¬R son literales.

**Forma normal disyuntiva.** Una fórmula F se dice que está en forma normal disyuntiva si y sólo si F es de la forma F = ($F_1 \lor F_2 \lor \ldots \lor F_n$), donde cada $F_i$ es una conjunción de literales.

**Forma normal conjuntiva.** Una fórmula F se dice que está en forma normal conjuntiva si y sólo si F es de la forma F = ($F_1 \land F_2 \land \ldots \land F_n$), donde cada $F_i$ es una disyunción de literales.

**Método de transformación**

**Paso 1**: Usar las leyes: $F \leftrightarrow G = (F \to G) \land (G \to F)$ y $F \to G = \neg F \lor G$ para eliminar las conectivas lógicas $\leftrightarrow$ y $\to$.

**Paso 2**: Usar repetidamente la ley $\neg(\neg F) = F$ y las leyes de De Morgan: $\neg(F \land H) = \neg F \lor \neg G$ y $\neg(F \lor H) = \neg F \land \neg G$ para disminuir el alcance de la negación a un único literal.

**Paso 3**: Usar de forma repetida las leyes distributivas:
$F \land (G \lor H) = (F \land G) \lor (F \land H)$ y $F \lor (G \land H) = (F \lor G) \land (F \lor H)$ y las otras leyes para obtener la forma normal. Es decir:

1. Cualquier función de conmutación de **n** variables **F(A, B, C,...)**, se puede expresar como una suma normal de productos utilizando los siguientes postulados:
   $A \land T = A$
   $A \lor \neg A = T$
   $A \land (B \lor \neg B) = (A \land B) \lor (A \land \neg B)$

2. Cualquier función de conmutación de **n** variables **F(A, B, C,...)**, se puede expresar como un producto normal de sumas, utilizando los siguientes postulados:
   $A \lor F = A$
   $A \land \neg A = F$
   $A \lor (B \land \neg B) = (A \lor B) \land (A \lor \neg B)$

**TÉRMINO PRODUCTO.** Conjunto de literales relacionadas por la conectiva $\land$

$P \land Q \land R, Q \land R \land S, P \land Q \land S$

**TÉRMINO SUMA.** Conjunto de literales relacionadas por la conectiva $\lor$

$P \lor Q \lor R, Q \lor R \lor S, P \lor Q \lor S$

**TÉRMINO NORMAL.** Un término producto o suma en el cual ninguna literal aparece más de una vez
- **Producto normal**
- **Suma normal**

**TÉRMINO CANÓNICO.** Término normal que contiene tantas literales como variables la función.

**Producto canónico o minitérmino.**
$P \land Q \land R, P \land Q \land \neg R, P \land \neg Q \land R$ ( para tres variables)

**Suma canónica o maxitérmino.**
$P \lor Q \lor R, P \lor Q \lor \neg R, P \lor \neg Q \lor R$ ( para tres variables)

**FORMA SUMA DE PRODUCTOS.** Una suma de términos producto (**MINITÉRMINO**) de una función.
$$F(A,B,C) = \sum\nolimits_{minit\acute{e}rminos}(\ ) = \sum\nolimits_m (\ )$$
F (P,Q,R)= (¬P $\land$ ¬Q $\land$ ¬R) $\lor$ (¬P $\land$ ¬Q $\land$ R) $\lor$ (¬P $\land$ Q $\land$ ¬R)

**FORMA PRODUCTO DE SUMAS.** Un producto de términos suma (**MAXITÉRMINOS**) de una función.
$$F(A,B,C) = \prod\nolimits_{MAXIT\acute{E}RMINOS}(\ ) = \prod\nolimits_M (\ )$$
F (P,Q,R)= (P $\lor$ Q $\lor$ R) $\land$ (P $\lor$ Q $\lor$ ¬R) $\land$ (P $\lor$ ¬Q $\lor$ R)

**FORMA CANÓNICA DE UNA FUNCIÓN.** Es aquella en que todos los términos son canónicos y aparecen una sola vez. Se tienen dos formas:
1. Suma de productos canónicos o suma de **MINITÉRMINOS**.
$$F(\ ) = \sum\nolimits_m (\ )$$
2. Producto de sumas canónicas o producto de **MAXITÉRMINOS**.
$$F(\ ) = \prod\nolimits_M (\ )$$

# Ejercicios de tautologías por medio de propiedades

**Ejercicio 1.**

**( P ∧ (P →Q))→Q ≡ T**
≡ (P ∧ (¬P ∨Q)) → Q        EL 1
≡ ¬ (P ∧ (¬P ∨ Q)) ∨ Q     EL 1
≡ (¬ P ∨ (P ∧ ¬Q)) ∨ Q     Ley de Morgan
≡ (¬ P ∨ Q) ∨ (P ∧ ¬Q)     Asociativa, Conmutativa
≡ ¬(P ∧ ¬Q) ∨ (P ∧ ¬Q)     Ley de Morgan
≡ T                         Negación

Otra forma
≡ (P ∧ (¬P ∨Q)) → Q             EL 1
≡ ¬ (P ∧ (¬P ∨ Q)) ∨ Q          EL 1
≡ (¬ P ∨ (P ∧ ¬Q)) ∨ Q          Ley de Morgan
≡ ((¬ P ∨ P) ∧ (¬ P∨ ¬Q)) ∨ Q   Distributiva
≡ ( T ∧ (¬ P∨ ¬Q)) ∨ Q          Negación
≡ (¬ P∨ ¬Q) ∨ Q                 Identidad
≡ ¬ P∨ (¬Q ∨ Q)                 Asociativa
≡ ¬ P∨ T                        Negación
≡ T                             Dominación

**Ejercicio 2.**

**((P →Q) ∧ ¬Q)→ ¬Q ≡ T**
≡ ¬ ((¬P ∨ Q) ∧ ¬Q) ∨ ¬Q       EL1
≡ (¬ (¬P ∨ Q) ∨ ¬ ¬Q) ∨ ¬Q     Ley de Morgan
≡ ((¬ ¬P ∧ ¬ Q) ∨ Q ) ∨ ¬Q     Ley de Morgan, Doble Negación
≡ ((P ∧ ¬ Q) ∨ Q ) ∨ ¬Q        Doble Negación
≡ ((P ∨ Q) ∧ (¬ Q ∨ Q )) ∨ ¬Q  Distributiva
≡ ((P ∨ Q) ∧ T ) ∨ ¬Q          Negación
≡ (P ∨ Q) ∨ ¬Q                 Identidad
≡ P ∨ (Q ∨ ¬Q)                 Asociativa
≡ P ∨ T                        Negación
≡ T                            Dominación

Otra forma

≡ ¬ ((¬P ∨ Q) ∧ ¬Q) ∨ ¬Q       EL1
≡ (¬ (¬P ∨ Q) ∨ ¬ ¬Q) ∨ ¬Q     Ley de Morgan
≡ ((¬ ¬P ∧ ¬ Q) ∨ Q ) ∨ ¬Q     Ley de Morgan, Doble Negación

≡ ((P ∧ ¬ Q) ∨ Q ) ∨ ¬Q          Doble Negación
≡ (P ∧ ¬Q) ∨ (Q ∨ ¬Q)            Asociativa
≡ (P ∧ ¬Q) ∨ T                   Negación
≡ T                              Dominación

**Ejercicio 3.**

**((P →Q) ∧ ¬(Q→R))→ (P→Q) ≡ T**
≡ ¬ ((¬P ∨ Q) ∧ ¬(¬Q∨R)) ∨ (¬P∨ Q)          EL1
≡ (¬ (¬P ∨ Q) ∨ ¬¬(¬Q∨R)) ∨ (¬P∨ Q)         Ley de Morgan
≡ ((¬ ¬P ∧ ¬Q) ∨ (¬Q∨R)) ∨ (¬P∨ Q)          Ley de Morgan, Doble Negación
≡ ((P ∧ ¬Q) ∨ (¬Q∨R)) ∨ (¬P∨ Q)             Doble Negación
≡ ((P ∧ ¬Q) ∨ ¬Q )∨ ( ¬P ∨ R ∨ Q)           Asociativa
≡ ¬Q ∨ ( ¬P ∨ R ∨ Q)                        Absorción
≡ (¬Q ∨ Q) ∨ (¬P ∨ R)                       Asociativa
≡ T ∨ (¬P ∨ R)                              Negación
≡ T                                         Dominación

Otra forma

≡ ¬ ((¬P ∨ Q) ∧ ¬(¬Q∨R)) ∨ (¬P∨ Q)          EL1
≡ ((P ∧ ¬Q) ∨ (¬Q∨R)) ∨ (¬P∨ Q)             Ley de Morgan
≡ ((P ∧ ¬Q) ∨ (¬P∨R) ) ∨ (¬Q∨ Q)            Asociativa
≡ ((P ∧ ¬Q) ∨ (¬P∨R) ) ∨ T                  Negación
≡ T                                         Dominación

Otra forma
≡ ¬ ((¬P ∨ Q) ∧ ¬(¬Q∨R)) ∨ (¬P∨ Q)          EL1
≡ ((P ∧ ¬Q) ∨ (¬Q∨R)) ∨ (¬P∨ Q)             Ley de Morgan
≡ ((P ∧ ¬Q) ∨ (¬P∨ Q)) ∨ (¬Q∨R)             Asociativa
≡ T ∨ (¬Q∨R)                                Negación
≡ T                                         Dominación

**Ejercicio 4.**

**P →(( P ∨ Q → R )→ (P →R))**

≡ ¬P ∨ (¬ ( ¬ (P ∨ Q) ∨ R ) ∨ (¬P ∨ R))     EL1
≡ ¬P ∨ (¬¬ (P ∨ Q) ∧ ¬ R ) ∨ (¬P ∨ R))      Ley de Morgan
≡ ¬P ∨ ((P ∨ Q) ∧ ¬ R ) ∨ (¬P ∨ R))         Doble Negación
≡ ((P ∨ Q) ∧ ¬ R ) ∨ (¬P ∨ ¬P ∨ R)          Asociativa
≡ ((P ∨ Q) ∧ ¬ R ) ∨ (¬P ∨ R)               Idem
≡ ((P ∧¬ R) ∨ (Q ∧ ¬ R )) ∨ (¬P ∨ R)        Distributiva
≡ ((P ∧¬ R) ∨ (¬P ∨ R))∨ (Q ∧ ¬ R )         Asociativa

| | |
|---|---|
| ≡ T ∨ (Q ∧ ¬R) | Negación |
| ≡ T | Dominación |

Otra forma

| | |
|---|---|
| ≡ ¬P ∨ (¬( ¬(P ∨ Q) ∨ R) ∨ (¬P ∨ R)) | EI1 |
| ≡ ¬P ∨ (((P ∨ Q) ∧ ¬R) ∨ (¬P ∨ R)) | Ley de Morgan |
| ≡ (¬P ∨ ¬P) ∨ (((P ∨ Q) ∧ ¬R) ∨ R) | Asociativa |
| ≡ ¬P ∨ (((P ∨ Q) ∨ R) ∧ (¬R ∨ R)) | Idem, Distributiva |
| ≡ ¬P ∨ (((P ∨ Q) ∨ R) ∧ T) | Negación |
| ≡ ¬P ∨ ((P ∨ Q) ∨ R) | Identidad |
| ≡ (¬P ∨ P) ∨ (Q ∨ R) | Asociativa |
| ≡ T ∨ (Q ∨ R) | Negación |
| ≡ T | Dominación |

**Ejercicio 5.**

**(P → R) ∧ (Q → S) → (P ∧ Q → R ∧ S)**

| | |
|---|---|
| ≡ ¬((¬P ∨ R) ∧ (¬Q ∨ S)) ∨ (¬(P ∧ Q) ∨ (R ∧ S)) | EI1 |
| ≡ ((P ∧ ¬R) ∨ (Q ∧ ¬S)) ∨ ((¬P ∨ ¬Q) ∨ (R ∧ S)) | Ley de Morgan |
| ≡ ((P ∧ ¬R) ∨ ¬P) ∨ ((Q ∧ ¬S) ∨ ¬Q) ∨ (R ∧ S) | Asociativa |
| ≡ ((P ∨ ¬P) ∧ (¬R ∨ ¬P)) ∨ ((Q ∨ ¬Q) ∧ (¬S ∨ ¬Q)) ∨ (R ∧ S) | Distributiva |
| ≡ (T ∧ (¬R ∨ ¬P)) ∨ (T ∧ (¬S ∨ ¬Q)) ∨ (R ∧ S) | Negación |
| ≡ (¬R ∨ ¬P) ∨ (¬S ∨ ¬Q) ∨ (R ∧ S) | Identidad |
| ≡ ((¬R ∨ ¬S) ∨ (R ∧ S)) ∨ (¬P ∨ ¬Q) | Asociativa |
| ≡ T ∨ (¬P ∨ ¬Q) | Negación |
| ≡ T | Dominación |

**Ejercicio 6.**

**(((P ∧ Q) → R) ∧ ((Q → R) → S) ∧ P) → S**

| | |
|---|---|
| ≡ ¬((¬(P ∧ Q) ∨ R) ∧ (¬(¬Q ∨ R) ∨ S) ∧ P) ∨ S | EI1 |
| ≡ (((P ∧ Q) ∧ ¬R) ∨ ((¬Q ∨ R) ∧ ¬S) ∨ ¬P) ∨ S | Ley de Morgan |
| ≡ (((P ∧ (Q ∧ ¬R)) ∨ ¬P) ∨ ((¬Q ∨ R) ∧ ¬S) ∨ S)) | Asociativa |
| ≡ ((P ∨ ¬P) ∧ ((Q ∧ ¬R) ∨ ¬P)) ∨ (((¬Q ∨ R) ∨ S) ∧ (¬S ∨ S)) | Distributiva |
| ≡ (T ∧ ((Q ∧ ¬R) ∨ ¬P)) ∨ (((¬Q ∨ R) ∨ S) ∧ T) | Negación |
| ≡ ((Q ∧ ¬R) ∨ ¬P)) ∨ ((¬Q ∨ R) ∨ S) | Identidad |
| ≡ ((Q ∧ ¬R) ∨ (¬Q ∨ R)) ∨ (¬P ∨ S) | Asociativa |
| ≡ T ∨ (¬P ∨ S) | Negación |
| ≡ T | Dominación |

# Anexo 5

*Material didáctico con videos*

## Ejercicios de Formas Normales por medio de propiedades.

1. **Expresar la función f (A, B, C) = A $\vee \overline{B} \wedge$ C en una suma de productos.**

$\equiv$ A$\vee \overline{B} \wedge$C

$\equiv$ A$\vee$ ( $\overline{B} \wedge$ C) Asociación  $\longleftarrow$  FND

$\equiv$ (A $\wedge$ T $\wedge$ T) $\vee$ (T $\wedge \overline{B} \wedge$C) identidad

$\equiv$ (A $\wedge$ (B$\vee\overline{B}$) $\wedge$ (C$\vee\overline{C}$) ) $\vee$ ( (A$\vee$A') $\wedge \overline{B} \wedge$C) negación

$\equiv$ (A$\wedge$B$\wedge$C)$\vee$(A$\wedge$B$\wedge\overline{C}$)$\vee$(A$\wedge\overline{B}\wedge$C)$\vee$(A$\wedge\overline{B}\wedge\overline{C}$)$\vee$(A$\wedge\overline{B}\wedge$C)$\vee$($\overline{A}\wedge\overline{B} \wedge$C) distributividad

$\equiv$ (A$\wedge$B$\wedge$C)$\vee$(A$\wedge$B$\wedge\overline{C}$)$\vee$ (A$\wedge\overline{B}\wedge$C) $\vee$(A$\wedge\overline{B}\wedge\overline{C}$)$\vee$($\overline{A}\wedge\overline{B} \wedge$C) idempotencia

$\equiv$ ($\overline{A}\wedge\overline{B} \wedge$C) $\vee$(A$\wedge\overline{B}\wedge\overline{C}$)$\vee$(A$\wedge\overline{B}\wedge$C)$\vee$ (A$\wedge$B$\wedge\overline{C}$)$\vee$(A$\wedge$B$\wedge$C) Conmutatividad  $\longleftarrow$  FNDP

 f (A, B ,C) =  $\sum$m =(1,4,5,6,7)

**Expresar la función F(A, B,C)= A$\vee \overline{B} \wedge$C en un producto de sumas.**

$\equiv$ A$\vee \overline{B} \wedge$C

$\equiv$ ( A$\vee \overline{B}$ ) $\wedge$ (A$\vee$ C) distributividad  $\longleftarrow$  FNC

$\equiv$ ( A$\vee \overline{B} \vee$ F ) $\wedge$ (A$\vee$ F $\vee$C) identidad

$\equiv$ ( A$\vee \overline{B} \vee$ (C$\wedge\overline{C}$) ) $\wedge$ (A$\vee$ (B$\wedge\overline{B}$) $\vee$C) negación

$\equiv$ ( A$\vee \overline{B} \vee$ C) $\wedge$( A$\vee \overline{B} \vee \overline{C}$) $\wedge$ (A$\vee$ B$\vee$C ) $\wedge$ (A$\vee \overline{B} \vee$C) distributividad

$\equiv$ ( A$\vee \overline{B} \vee$ C) $\wedge$( A$\vee \overline{B} \vee \overline{C}$) $\wedge$ (A$\vee$ B$\vee$C ) idempotencia

$\equiv$ (A$\vee$ B$\vee$C ) $\wedge$( A$\vee \overline{B} \vee$ C) $\wedge$ ( A$\vee \overline{B} \vee \overline{C}$)  Conmutatividad  $\longleftarrow$  FNCP

   f (A, B, C) = $\Pi_M$ = (0,2,3)

$\therefore$ **Contingencia**

| Puedes ver la explicación en video: | 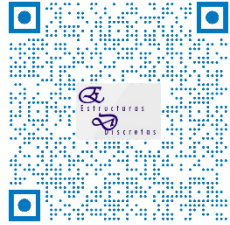 | 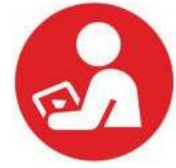 |
|---|---|---|

2. **Obteniendo su FNCP de la siguiente formula $(((P \wedge Q) \rightarrow P) \rightarrow (Q \vee R)) \wedge (\overline{Q} \wedge \overline{R})$**

$\equiv (((P \wedge Q) \rightarrow P) \rightarrow (Q \vee R)) \wedge (\overline{Q} \wedge \overline{R})$

$\equiv (\overline{(\overline{(P \wedge Q)}P)} \vee (Q \vee R)) \wedge (\overline{Q} \wedge \overline{R}))$ *EL 1*

$\equiv (((P \wedge Q)\overline{P}) \vee (Q \vee R)) \wedge (\overline{Q} \wedge \overline{R})$ Ley de Morgan

$\equiv (F \vee (Q \vee R)) \wedge (\overline{Q} \wedge \overline{R})$ Conmutatividad, Asociación, Negación

$\equiv (Q \vee R) \wedge (\overline{Q} \wedge \overline{R})$ identidad

$\equiv (Q \vee R) \wedge \overline{Q} \wedge \overline{R}$ Asociación   ⟵   **FNC**

$\equiv (F \vee Q \vee R) \wedge (F \vee \overline{Q} \vee F) \wedge (F \vee F \vee \overline{R})$ identidad

$\equiv ((P \wedge \overline{P}) \vee Q \vee R) \wedge ((P \wedge \overline{P}) \vee \overline{Q} \vee (R \wedge \overline{R})) \wedge ((P \wedge \overline{P}) \vee (Q \wedge \overline{Q}) \vee \overline{R})$ negación

$\equiv (P \vee Q \vee R) \wedge (\overline{P} \vee Q \vee R) \wedge (P \vee \overline{Q} \vee R) \wedge (P \vee \overline{Q} \vee \overline{R}) \wedge (\overline{P} \vee \overline{Q} \vee R) \wedge (\overline{P} \vee \overline{Q} \vee \overline{R}) \wedge (P \vee Q \vee \overline{R}) \wedge (P \vee \overline{Q} \vee \overline{R}) \wedge (\overline{P} \vee Q \vee \overline{R}) \wedge (\overline{P} \vee \overline{Q} \vee \overline{R})$ distributividad

$\equiv (P \vee Q \vee R) \wedge (\overline{P} \vee Q \vee R) \wedge (P \vee \overline{Q} \vee R) \wedge (P \vee \overline{Q} \vee \overline{R}) \wedge (\overline{P} \vee Q \vee R) \wedge (\overline{P} \vee \overline{Q} \vee R) \wedge (P \vee Q \vee \overline{R}) \wedge (\overline{P} \vee Q \vee \overline{R})$ idempotencia

$\equiv (P \vee Q \vee R) \wedge (P \vee Q \vee \overline{R}) \wedge (P \vee \overline{Q} \vee R) \wedge (P \vee \overline{Q} \vee \overline{R}) \wedge (\overline{P} \vee Q \vee R) \wedge (\overline{P} \vee Q \vee \overline{R}) \wedge (\overline{P} \vee \overline{Q} \vee R) \wedge (\overline{P} \vee \overline{Q} \vee \overline{R})$ conmutatividad   ⟵   FNCP

$f(P,Q,R) = \Pi_M = (0,1,2,3,4,5,6,7)$

∴ **Contradicción**

| Puedes ver la explicación en video: | 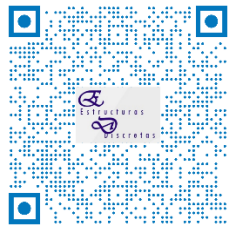 | 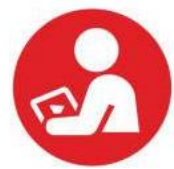 |
|---|---|---|

**3. Verificar por medio de formas normales si son inconsistencia, tautologías o contradicciones**   $p \land (\overline{q} \lor (r \land (\overline{s} \land ( r \lor ( \overline{q} \lor p)))))$

$\equiv p \land (\overline{q} \lor (r \land (\overline{s} \land ( r \lor ( \overline{q} \lor p)))))$

$\equiv p \land (\overline{q} \lor (\overline{s} \land (r \land ( r \lor p \lor \overline{q}) )))$ Asociatividad, conmutatividad

$\equiv p \land (\overline{q} \lor (r \land \overline{s}))$ Absorción, conmutatividad

$\equiv (p \land \overline{q}) \lor (p \land r \land \overline{s})$ Distributividad, asociación  ⬅ **FND**

$\equiv (p \land \overline{q} \land T \land T) \lor (p \land T \land r \land \overline{s})$ Identidad

$\equiv (p \land \overline{q} \land (r \lor \overline{r}) \land (s \lor \overline{s})) \lor (p \land (q \lor \overline{q}) \land r \land \overline{s})$ Negación

$\equiv (p \land \overline{q} \land r \land s) \lor (p \land \overline{q} \land r \land \overline{s}) \lor (p \land \overline{q} \land \overline{r} \land s) \lor (p \land \overline{q} \land \overline{r} \land \overline{s}) \lor (p \land q \land r \land \overline{s}) \lor (p \land \overline{q} \land r \land \overline{s})$ Distributividad

$\equiv (p \land \overline{q} \land r \land s) \lor (p \land \overline{q} \land r \land \overline{s}) \lor (p \land \overline{q} \land \overline{r} \land s) \lor (p \land \overline{q} \land \overline{r} \land \overline{s}) \lor (p \land q \land r \land \overline{s})$ Idempotencia

$\equiv (p \land \overline{q} \land \overline{r} \land \overline{s}) \lor (p \land \overline{q} \land \overline{r} \land s) \lor (p \land \overline{q} \land r \land \overline{s}) \lor (p \land \overline{q} \land r \land s) \lor (p \land q \land r \land \overline{s})$ Conmutatividad  ⬅ **FNDP**

**f(A, B,C) =** ∑m =(8,9,10,11,14)

De este paso podemos obtener su forma normal conjuntiva

$\equiv p \land (\overline{q} \lor r) \land (\overline{q} \lor \overline{s})$ Distributividad, asociativida  ⬅ **FNC**

$\equiv (p \lor F \lor F \lor F) \land (F \lor \overline{q} \lor r \lor F) \land (F \lor \overline{q} \lor F \lor \overline{s})$ identidad

$\equiv (p \lor (q \land \overline{q}) \lor (r \land \overline{r}) \lor (s \land \overline{s})) \land ((p \land \overline{p}) \lor \overline{q} \lor r \lor (s \land \overline{s})) \land ((p \land \overline{p}) \lor \overline{q} \lor (r \land \overline{r}) \lor \overline{s})$ Negación

$\equiv (p \lor q \lor r \lor s) \land (p \lor q \lor r \lor \overline{s}) \land (p \lor \overline{q} \lor r \lor s) \land (p \lor \overline{q} \lor r \lor \overline{s}) \land (p \lor q \lor \overline{r} \lor s) \land (p \lor q \lor \overline{r} \lor \overline{s}) \land (p \lor \overline{q} \lor \overline{r} \lor s) \land (p \lor \overline{q} \lor \overline{r} \lor \overline{s}) \land (p \lor \overline{q} \lor r \lor s) \land (p \lor \overline{q} \lor r \lor \overline{s}) \land (\overline{p} \lor \overline{q} \lor r \lor s) \land (\overline{p} \lor \overline{q} \lor r \lor \overline{s}) \land (p \lor \overline{q} \lor r \lor \overline{s}) \land (p \lor \overline{q} \lor \overline{r} \lor \overline{s}) \land (\overline{p} \lor \overline{q} \lor r \lor \overline{s}) \land (\overline{p} \lor \overline{q} \lor \overline{r} \lor \overline{s})$ Distributividad

$\equiv (p \lor q \lor r \lor s) \land (p \lor q \lor r \lor \overline{s}) \land (p \lor \overline{q} \lor r \lor s) \land (p \lor \overline{q} \lor r \lor \overline{s}) \land (p \lor q \lor \overline{r} \lor s) \land (p \lor q \lor \overline{r} \lor \overline{s}) \land (p \lor \overline{q} \lor \overline{r} \lor s) \land (p \lor \overline{q} \lor \overline{r} \lor \overline{s}) \land (\overline{p} \lor \overline{q} \lor r \lor s) \land (\overline{p} \lor \overline{q} \lor r \lor \overline{s}) \land (\overline{p} \lor \overline{q} \lor \overline{r} \lor \overline{s})$ Idempotencia

≡ (p∨q∨r ∨s) ∧ (p∨q∨r∨ $\bar{s}$) ∧(p∨q∨$\bar{r}$ ∨s) ∧ (p∨q∨$\bar{r}$∨ $\bar{s}$) ∧ (p∨$\bar{q}$∨r ∨s) ∧ (p∨$\bar{q}$∨r ∨$\bar{s}$) ∧(p∨$\bar{q}$∨$\bar{r}$ ∨s)∧ (p∨$\bar{q}$∨$\bar{r}$ ∨$\bar{s}$) ∧($\bar{p}$∨$\bar{q}$∨r∨s) ∧($\bar{p}$∨$\bar{q}$∨r∨$\bar{s}$) ∧ ($\bar{p}$∨$\bar{q}$ ∨$\bar{r}$∨ $\bar{s}$) Conmutatividad ⬅

**FNCP**

**f(p,q,r) = $\Pi_M$ = (0,1,2,3,4,5,6,7,12,13,15)**

∴ **Contingencia**

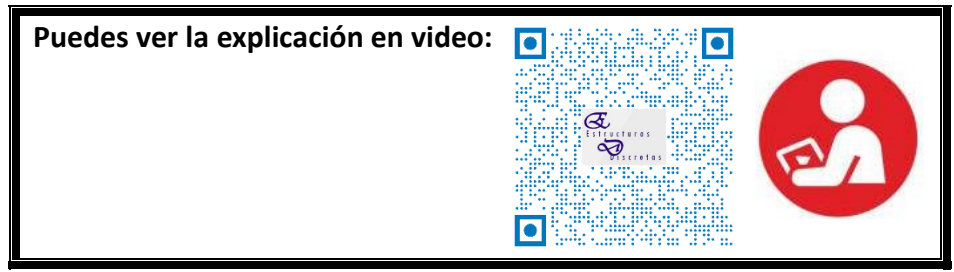
Puedes ver la explicación en video:

**4. Demostrar que la siguiente formula es una inconsistencia f(p,q,r) = ($\bar{p}$ ∨$\bar{q}$) ∧ (p∨ r)**

Obteniendo su FNCP

≡ ($\bar{p}$ ∨$\bar{q}$) ∧ (p∨ r)

≡ (($\bar{p}$ ∨$\bar{q}$) ∨ F) ∧ (p∨ F ∨ r) identidad

≡ (($\bar{p}$ ∨$\bar{q}$) ∨ (r∧$\bar{r}$)) ∧ ((p∨ (q∧$\bar{q}$) ∨ r)) negación

≡ ( ($\bar{p}$ ∨$\bar{q}$ ∨ r) ∧ ($\bar{p}$ ∨$\bar{q}$ ∨$\bar{r}$) ) ∧ ( (p∨q∨ r) ∧ (p∨$\bar{q}$ ∨ r) ) distributividad

≡ ($\bar{p}$ ∨$\bar{q}$ ∨ r) ∧ ($\bar{p}$ ∨$\bar{q}$ ∨$\bar{r}$) ∧ (p∨q∨ r) ∧ (p∨$\bar{q}$ ∨ r) Asociación

≡ (p∨q∨ r) ∧ (p∨$\bar{q}$ ∨ r) ∧($\bar{p}$ ∨$\bar{q}$ ∨ r) ∧ ($\bar{p}$ ∨$\bar{q}$ ∨$\bar{r}$) Conmutatividad ⬅ FNCP (Forma Normal Conjuntiva Principal)

**f(p,q,r) = $\Pi_M$ = (0,2,6,7)**

Obteniendo su FNDP

f(p,q,r) = ($\bar{p}$ ∨$\bar{q}$) ∧ (p∨ r)

≡ ($\bar{p}$ ∨$\bar{q}$) ∧ (p∨ r)

≡ (($\bar{p}$ ∨$\bar{q}$) ∧ p)∨ (($\bar{p}$ ∨$\bar{q}$) ∧ r) distributividad

≡ ( ($\bar{p}$ ∧ p) ∨ ($\bar{q}$ ∧ p)) ∨ ( ($\bar{p}$ ∧ r )∨($\bar{q}$ ∧ r) ) distributividad

≡ ( F ∨ ($\overline{q}$ ∧ p)) ∨(($\overline{p}$ ∧ r )∨($\overline{q}$ ∧ r ) )negación

≡ ($\overline{q}$ ∧ p) ∨($\overline{p}$ ∧ r )∨($\overline{q}$ ∧ r) identidad, Asociación

≡ (p∧$\overline{q}$ ∧ T) ∨($\overline{p}$∧ T∧ r )∨(T ∧$\overline{q}$ ∧ r) identidad

≡ (p∧$\overline{q}$ ∧ (r ∨$\overline{r}$)) ∨($\overline{p}$∧ (q ∨$\overline{q}$ ) ∧ r )∨( (p∨$\overline{p}$) ∧ $\overline{q}$ ∧ r) negación

≡ ( (p∧$\overline{q}$ ∧ r) ∨ (p∧$\overline{q}$ ∧ $\overline{r}$)) ∨( ($\overline{p}$∧ q∧ r) ∨($\overline{p}$∧$\overline{q}$ ∧ r )) ∨ ( (p∧$\overline{q}$ ∧ r)∨($\overline{p}$∧$\overline{q}$ ∧ r) ) distributividad

≡ (p∧$\overline{q}$ ∧ r) ∨ (p∧$\overline{q}$ ∧ $\overline{r}$) ∨($\overline{p}$∧ q∧ r) ∨ ($\overline{p}$∧$\overline{q}$ ∧ r ) Asociación

≡ ($\overline{p}$∧$\overline{q}$ ∧ r) ∨ ($\overline{p}$∧ q∧ r) ∨ (p∧$\overline{q}$ ∧ $\overline{r}$) ∨  (p∧$\overline{q}$ ∧ r) Conmutatividad ⬅ FNDP (Forma Normal Disyuntiva Principal)

**f(p,q,r) = Σ = (1,3,4,5)**

∴ **Contingencia**

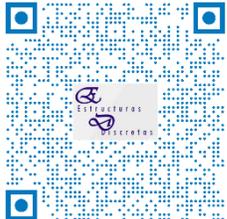

Puedes ver la explicación en video:

5. Obteniendo la FNCP     f(a,b,c,d) = (a∨$\overline{b}$∨ d) ∧ ((b∧d)∨$\overline{a∧ \overline{c}}$)

≡ (a∨$\overline{b}$∨ d) ∧ ((b∧d)∨ $\overline{a∧ \overline{c}}$)

≡ (a∨$\overline{b}$∨ d) ∧ ((b∧d)∨ ($\overline{a}$∨ c)) Morgan

≡ (a∨$\overline{b}$∨ d) ∧ ((b ∨ ($\overline{a}$∨ c)) ∧ (d ∨ ($\overline{a}$∨ c)) Distributividad

≡ (a∨$\overline{b}$∨ d) ∧ (b ∨ $\overline{a}$∨ c) ∧ (d ∨ $\overline{a}$∨ c) Asociatividad

≡ (a∨$\overline{b}$∨ d) ∧ ($\overline{a}$∨b ∨ c) ∧ ($\overline{a}$∨ c∨ d) Conmutatividad

≡ (a∨$\overline{b}$∨ F ∨d) ∧ ($\overline{a}$∨b ∨ c ∨ F) ∧ ($\overline{a}$∨ F ∨ c∨ d) identidad

≡ (a∨$\overline{b}$∨ (c∧ $\overline{c}$) ∨d) ∧ ($\overline{a}$∨b ∨ c ∨ (d∧ $\overline{d}$) ) ∧ ($\overline{a}$∨ (b∧$\overline{b}$) ∨ c∨ d) Negación

$\equiv (a \vee \overline{b} \vee c \vee d) \wedge (a \vee \overline{b} \vee \overline{c} \vee d) \wedge (\overline{a} \vee b \vee c \vee d) \wedge (\overline{a} \vee b \vee c \vee \overline{d}) \wedge (\overline{a} \vee b \vee c \vee d) \wedge (\overline{a} \vee \overline{b} \vee c \vee d)$
Distributividad

$\equiv (a \vee \overline{b} \vee c \vee d) \wedge (a \vee \overline{b} \vee \overline{c} \vee d) \wedge (\overline{a} \vee b \vee c \vee d) \wedge (\overline{a} \vee b \vee c \vee \overline{d}) \wedge (\overline{a} \vee \overline{b} \vee c \vee d)$ Idempotencia

$f(a, b, c, d) = (a \vee \overline{b} \vee c \vee d) \wedge (a \vee \overline{b} \vee \overline{c} \vee d) \wedge (\overline{a} \vee b \vee c \vee d) \wedge (\overline{a} \vee b \vee c \vee \overline{d}) \wedge (\overline{a} \vee \overline{b} \vee c \vee d)$
Conmutatividad → **FNCP**

$f(a, b, c, d) = \Pi_M (4, 6, 8, 9, 12)$

Obteniendo la FNDP

$\equiv (a \vee \overline{b} \vee d) \wedge ((b \wedge d) \vee \overline{\overline{a} \wedge \overline{c}})$

$\equiv (a \vee \overline{b} \vee d) \wedge ((b \wedge d) \vee (\overline{a} \vee c))$ Morgan

$\equiv ((a \vee \overline{b} \vee d) \wedge (b \wedge d)) \vee ((a \vee \overline{b} \vee d) \wedge (\overline{a} \vee c))$ Distributividad

$\equiv ((a \vee \overline{b} \vee d) \wedge b \wedge d) \vee ((a \vee \overline{b} \vee d) \wedge (\overline{a} \vee c))$ Asociatividad

$\equiv ((a \vee \overline{b} \vee d) \wedge d \wedge b) \vee ((a \vee \overline{b} \vee d) \wedge (\overline{a} \vee c))$ Conmutatividad

$\equiv (d \wedge b) \vee ((a \vee \overline{b} \vee d) \wedge (\overline{a} \vee c))$ Absorción

$\equiv (d \wedge b) \vee (((a \vee \overline{b} \vee d) \wedge \overline{a}) \vee ((a \vee \overline{b} \vee d) \wedge c))$ Distributividad

$\equiv (d \wedge b) \vee ((a \wedge \overline{a}) \vee (\overline{b} \wedge \overline{a}) \vee (d \wedge \overline{a})) \vee ((a \wedge c) \vee (\overline{b} \wedge c) \vee (d \wedge c))$ Distributividad

$\equiv (d \wedge b) \vee (a \wedge \overline{a}) \vee (\overline{b} \wedge \overline{a}) \vee (d \wedge \overline{a}) \vee (a \wedge c) \vee (\overline{b} \wedge c) \vee (d \wedge c)$ Asociatividad

$\equiv (d \wedge b) \vee F \vee (\overline{b} \wedge \overline{a}) \vee (d \wedge \overline{a}) \vee (a \wedge c) \vee (\overline{b} \wedge c) \vee (d \wedge c)$ Negación

$\equiv (d \wedge b) \vee (\overline{b} \wedge \overline{a}) \vee (d \wedge \overline{a}) \vee (a \wedge c) \vee (\overline{b} \wedge c) \vee (d \wedge c)$ identidad

$\equiv (b \wedge d) \vee (\overline{a} \wedge \overline{b}) \vee (\overline{a} \wedge d) \vee (a \wedge c) \vee (\overline{b} \wedge c) \vee (c \wedge d)$ Conmutatividad

$\equiv (T \wedge b \wedge T \wedge d) \vee (\overline{a} \wedge \overline{b} \wedge T \wedge T) \vee (\overline{a} \wedge T \wedge T \wedge d) \vee (a \wedge T \wedge c \wedge T) \vee (T \wedge \overline{b} \wedge c \wedge T) \vee (T \wedge T \wedge c \wedge d)$ identidad

$\equiv ((a \vee \overline{a}) \wedge b \wedge (c \vee \overline{c}) \wedge d) \vee (\overline{a} \wedge \overline{b} \wedge (c \vee \overline{c}) \wedge (d \vee \overline{d})) \vee (\overline{a} \wedge (b \vee \overline{b}) \wedge (c \vee \overline{c}) \wedge d) \vee (a \wedge (b \vee \overline{b}) \wedge c \wedge (d \vee \overline{d}))$
$\vee ((a \vee \overline{a}) \wedge \overline{b} \wedge c \wedge (d \vee \overline{d})) \vee ((a \vee \overline{a}) \wedge (b \vee \overline{b}) \wedge c \wedge d)$ Negación

$\equiv (a \wedge b \wedge c \wedge d) \vee (a \wedge b \wedge \overline{c} \wedge d) \vee (\overline{a} \wedge b \wedge c \wedge d) \vee (\overline{a} \wedge b \wedge \overline{c} \wedge d) \vee (\overline{a} \wedge \overline{b} \wedge c \wedge d) \vee (\overline{a} \wedge \overline{b} \wedge c \wedge \overline{d}) \vee (\overline{a} \wedge \overline{b} \wedge \overline{c} \wedge d)$
$\vee (\overline{a} \wedge \overline{b} \wedge \overline{c} \wedge \overline{d}) \vee (\overline{a} \wedge b \wedge c \wedge d) \vee (\overline{a} \wedge b \wedge \overline{c} \wedge d) \vee (\overline{a} \wedge \overline{b} \wedge c \wedge d) \vee (\overline{a} \wedge \overline{b} \wedge \overline{c} \wedge d) \vee (a \wedge b \wedge c \wedge d) \vee (a \wedge b \wedge c \wedge \overline{d})$
$\vee (a \wedge \overline{b} \wedge c \wedge d) \vee (a \wedge \overline{b} \wedge c \wedge \overline{d}) \vee (a \wedge \overline{b} \wedge c \wedge d) \vee (a \wedge \overline{b} \wedge c \wedge \overline{d}) \vee (\overline{a} \wedge \overline{b} \wedge c \wedge d) \vee (\overline{a} \wedge \overline{b} \wedge c \wedge \overline{d}) \vee (a \wedge b \wedge c \wedge d) \vee$
$(a \wedge \overline{b} \wedge c \wedge d) \vee (\overline{a} \wedge b \wedge c \wedge d) \vee (\overline{a} \wedge \overline{b} \wedge c \wedge d)$ Distributividad

≡ (a∧b∧c∧d)∨(a∧b∧$\bar{c}$∧d)∨($\bar{a}$∧b∧c∧d)∨($\bar{a}$∧b∧$\bar{c}$∧d)∨($\bar{a}$∧$\bar{b}$∧c∧d)∨($\bar{a}$∧$\bar{b}$∧c∧$\bar{d}$)∨($\bar{a}$∧$\bar{b}$∧$\bar{c}$∧d) ∨($\bar{a}$∧$\bar{b}$∧$\bar{c}$∧$\bar{d}$)∨(a∧b∧c∧$\bar{d}$) ∨(a∧$\bar{b}$∧c∧d)∨(a∧$\bar{b}$∧c∧$\bar{d}$) idempotencia

≡ ($\bar{a}$∧$\bar{b}$∧$\bar{c}$∧$\bar{d}$) ∨ ($\bar{a}$∧$\bar{b}$∧$\bar{c}$∧d) ∨ ($\bar{a}$∧$\bar{b}$∧c∧$\bar{d}$) ∨($\bar{a}$∧$\bar{b}$∧c∧d) ∨ ($\bar{a}$∧b∧$\bar{c}$∧d) ∨ ($\bar{a}$∧b∧c∧d) ∨ (a∧$\bar{b}$∧c∧$\bar{d}$) ∨(a∧$\bar{b}$∧c∧d) ∨ (a∧b∧$\bar{c}$∧d) ∨ (a∧b∧c∧$\bar{d}$) ∨ (a∧b∧c∧d) Conmutatividad → **FNDP**

f(a,b,c,d)=($\bar{a}$∧$\bar{b}$∧$\bar{c}$∧$\bar{d}$)∨($\bar{a}$∧$\bar{b}$∧$\bar{c}$∧d)∨($\bar{a}$∧$\bar{b}$∧c∧$\bar{d}$)∨($\bar{a}$∧$\bar{b}$∧c∧d)∨($\bar{a}$∧b∧$\bar{c}$∧d)∨($\bar{a}$∧b∧c∧d)∨(a∧$\bar{b}$∧c∧$\bar{d}$)∨(a∧$\bar{b}$∧c∧d)∨(a∧b∧$\bar{c}$∧d)∨(a∧b∧c∧$\bar{d}$)∨(a∧b∧c∧d) = $\Sigma_m$ **(0,1,2,3,5,7,10,11,13,14,15)**

∴ **Contingencia**

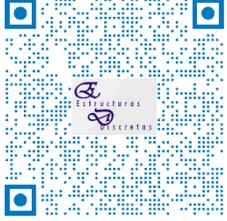

*Material con videos y Jamboard*

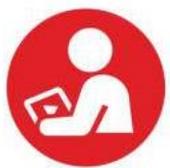

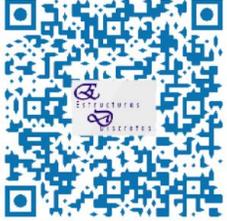

**P, P→(Q ∨ R), (Q ∨ R) →S ⇒ S**

$P \wedge (P \to (Q \vee R)) \wedge ((Q \vee R) \to S) \Rightarrow S$

Directo

Se parte $P \wedge (P \to (Q \vee R)) \wedge ((Q \vee R) \to S) \equiv T \ldots$ ①

De ①    $P \equiv T, \ldots$ ②

$P \to (Q \vee R) \equiv T$   ③

$(Q \vee R) \to S \equiv T$   ④

Sust ② en ③   $T \to (Q \vee R) \equiv T$

$(Q \vee R) \equiv T \ldots$ ⑤

Sust ⑤ en ④   $T \to S \equiv T$

$S \equiv T \ldots$ ⑥

Demostrando $Q \equiv S \equiv T$, sust ⑥

$T \equiv T$    ∴ CL válida

---

**P, P→(Q ∨ R), (Q ∨ R) →S ⇒ S**

Indirecto

Se tiene $S \equiv F \ldots$ ①

Demostrando $\overline{P \wedge (P \to (Q \vee R)) \wedge ((Q \vee R) \to S)} \equiv F$

$\equiv \overline{P \wedge (P \to (Q \vee R)) \wedge ((Q \vee R) \to S)}$

$\equiv P \wedge (\bar{P} \vee (Q \vee R)) \wedge (\overline{(Q \vee R)} \vee F)$   ELI, Sust ①

$\equiv P \wedge (\bar{P} \vee (Q \vee R)) \wedge \overline{(Q \vee R)}$   Ident

$\equiv (P \wedge \overline{(Q \vee R)}) \wedge (\bar{P} \vee (Q \vee R))$   Conm, Asoc

$\equiv \overline{(\bar{P} \vee (Q \vee R))} \wedge (\bar{P} \vee (Q \vee R))$   Morgan

$\equiv F$   Dom    ∴ CL válida

**P ∨ Q, P→R, Q→R ⇒ R**

$(P \vee Q) \wedge (P \to R) \wedge (Q \to R) \to R$

**Directo**

Se parte $(P \vee Q) \wedge (P \to R) \wedge (Q \to R) \equiv T \ldots ①$

De ① 
- $P \vee Q \equiv T \ldots ②$
- $P \to R \equiv T \ldots ③$
- $Q \to R \equiv T \ldots ④$

Trabajando con ② y ③

| | |
|---|---|
| $(P \vee Q) \wedge (P \to R) \equiv T$ | |
| $(P \vee Q) \wedge (\bar{P} \vee R) \equiv T$ | EL1 |
| $((P \vee Q) \wedge \bar{P}) \vee ((P \vee Q) \wedge R) \equiv T$ | Dist |
| $((P \wedge \bar{P}) \vee (Q \wedge \bar{P})) \vee (T \wedge R) \equiv T$ | Dist, Sust ② |
| $(F \vee (Q \wedge \bar{P})) \vee R \equiv T$ | Neg, ident |
| $(Q \wedge \bar{P}) \vee R \equiv T$ | ident |
| $(Q \vee R) \wedge (\bar{P} \vee R) \equiv T$ | Dist |
| $(Q \vee R) \wedge (P \to R) \equiv T$ | EL 1 |
| $(Q \vee R) \wedge T \equiv T$ | Sust ③ |
| $Q \vee R \equiv T \ldots ⑤$ | ident |

Trabajando con ④ y ⑤ 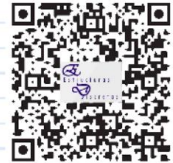

| | |
|---|---|
| $(Q \to R) \wedge (Q \vee R) \equiv T$ | |
| $(\bar{Q} \vee R) \wedge (Q \vee R) \equiv T$ | EL1 |
| $(\bar{Q} \wedge Q) \vee R \equiv T$ | Dist |
| $F \vee R \equiv T$ | Neg |
| $R \equiv T \ldots ⑥$ | ident |

Demostrando $①\equiv T$, sust ⑥

$R \equiv T$
$T \equiv T$

∴ CL válida

---

**P ∨ Q, P→R, Q→R ⇒ R**

**Indirecto**

Se tiene $Q \equiv R \equiv F \ldots ①$

Demostrando $(P \vee Q) \wedge (P \to R) \wedge (Q \to R) \equiv F$

| | |
|---|---|
| $\equiv (P \vee Q) \wedge (P \to R) \wedge (Q \to R)$ | |
| $\equiv (P \vee Q) \wedge (\bar{P} \vee F) \wedge (\bar{Q} \vee F)$ | EL1, Sust ① |
| $\equiv (P \vee Q) \wedge \bar{P} \wedge \bar{Q}$ | ident |
| $\equiv \overline{(P \vee Q)} \wedge (\bar{P} \wedge \bar{Q})$ | Morgan, Asoc |
| $\equiv F$ | Neg |

∴ CL válida

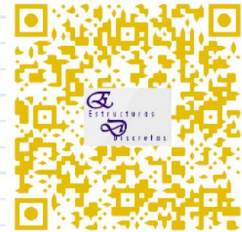

**P→Q, Q→R, ¬R ⇒ ¬P**

Definición
$(P \to Q) \land (Q \to R) \land \bar{R} \to \bar{P} \equiv T$
$\equiv \overline{(P \to Q) \land (Q \to R) \land \bar{R}} \lor \bar{P}$
$\equiv \overline{(\bar{P} \lor Q) \land (\bar{Q} \lor R) \land \bar{R}} \lor \bar{P}$  EL 1
$\equiv \overline{(\bar{P} \lor Q) \land ((\bar{Q} \land \bar{R}) \lor (R \land \bar{R}))} \lor \bar{P}$  Dist
$\equiv \overline{(\bar{P} \lor Q) \land ((\bar{Q} \land \bar{R}) \lor F)} \lor \bar{P}$  Neg
$\equiv \overline{(\bar{P} \lor Q) \land (\bar{Q} \land \bar{R})} \lor \bar{P}$  ident
$\equiv \overline{((\bar{P} \lor Q) \land \bar{Q}) \land \bar{R}} \lor \bar{P}$  Asoc
$\equiv \overline{((\bar{P} \land \bar{Q}) \lor (Q \land \bar{Q})) \land \bar{R}} \lor \bar{P}$  Dist
$\equiv \overline{(\bar{P} \land \bar{Q}) \land \bar{R}} \lor \bar{P}$  Neg, ident
$\equiv (P \lor Q) \lor R \lor \bar{P}$
$\equiv (P \lor \bar{P}) \lor Q \lor R$  Asoc, Conm
$\equiv T \lor Q \lor R$  Neg
$\equiv T$  Dom

∴ CL válida

Puedes ver la explicación en video:

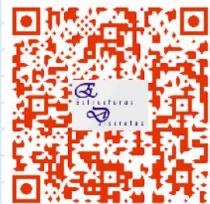

---

**P→Q, Q→R, ¬R ⇒ ¬P**

Directo
Se parte $(P \to Q) \land (Q \to R) \land \bar{R} \equiv T \ldots 1$
De ① $P \to Q \equiv T \ldots$ ②
$Q \to R \equiv T \ldots$ ③
$\bar{R} \equiv T \ldots$ ④
Sust ④ en ③   $Q \equiv F \ldots$ ⑤
Sust ⑤ en ②   $P \equiv F \ldots$ ⑥
Demostrando $Q \equiv \bar{P} \equiv T$, sust ⑥
$T \equiv \bar{F}$
$T \equiv T$   ∴ CL válido

Indirecto
Partir $Q \equiv F \equiv \bar{P} \ldots$ ①
Demostrando $(P \to Q) \land (Q \to R) \land \bar{R} \equiv F$
$\equiv (P \to Q) \land (Q \to R) \land \bar{R}$
$\equiv (\bar{P} \lor Q) \land (\bar{Q} \lor R) \land \bar{R}$  EL 1
$\equiv (F \lor Q) \land (\bar{Q} \lor R) \land \bar{R}$  Sust ①
$\equiv Q \land (\bar{Q} \lor R) \land \bar{R}$  ident
$\equiv (Q \land \bar{R}) \land (\bar{Q} \lor R)$  Asoc, Conm
$\equiv F$  Neg

∴ CL válida

# Anexo 6

*Diseño de aula virtual de Google Classroom*